\documentclass{article}

\usepackage{amssymb} 
\usepackage{amsmath} 
\usepackage{latexsym} 
\usepackage{theorem} 
\usepackage{color}
\usepackage{epsfig}

\theorembodyfont{\itshape} 

\newtheorem{theorem}{Theorem}[section]
\newtheorem{lemma}[theorem]{Lemma}

\newtheorem{corollary}[theorem]{Corollary}
\newtheorem{proposition}[theorem]{Proposition}

\newtheorem{definition}[theorem]{Definition}
\newtheorem{example}{Example}[section]

\newtheorem{problem}[theorem]{Problem}
\newtheorem{remark}[theorem]{Remark}
\newenvironment{proof}{{\par\addvspace{0.1cm}\noindent \bf Proof. }}{\hfill$\Box$\par\medskip} 
 
\def\n{m}

\def\an{{\alpha-\n}}
\def\a{\alpha}
\def\e{\varepsilon}
\def\Om{\varOmega}
\def\RR{\mathbb{R}}
\def\vect#1{\mbox{\boldmath $#1$}} 
\def\interior#1{\stackrel{\circ}{#1}}
\def\diam{{d}}
\def\rhosub#1{\mbox{\large $\rho $}_{\mbox{\small $\!{}_{#1}$}}}

\newcommand{\ti}{\;\;\makebox[0pt]{$\top$}\makebox[0pt]{\small $\cap$}\;\;}

\title{\bf Renormalization of 
potentials and generalized centers}
\author{Jun O'Hara}

\numberwithin{equation}{section}

\begin{document}

\maketitle

\begin{abstract} 
We generalize the Riesz potential of a compact domain in $\mathbb{R}^{m}$ by introducing a renormalization of the $r^{\alpha-m}$-potential for $\alpha\le0$. This can be considered as generalization of the dual mixed volumes of convex bodies as introduced by Lutwak. 
We then study the points where the extreme values of the (renormalized) potentials are attained. These points can be considered as a generalization of the center of mass.  
We also show that only balls give extreme values among bodied with the same volume. 
\end{abstract}




\section{Introduction}
Let $\Om$ be a compact set in $\RR^\n$ $(\n\ge2)$ which can be obtained as a closure of an open set. 
Consider a potential of the form 
\[V_\a(x)=\int_\Om{|x-y|}^{\an}\,d\mu(y),\]
where $\mu$ is the standard Lesbegue measure of $\RR^\n$. 
When $\a<\n$ and $x\in\Om$ it is a singular integral that is well-defined if $\a>0$. 
When $0<\a<m$ it is the {\em Riesz potential} of the characteristic function $\chi_\Om$ of $\Om$. 
In particular, it is (a constant times) the {\em Newton potential} when $\a=2$ and $\n\ge3$. 

We can apply the same renormalization process which was used to define energy functionals of knots (\cite{O1,O2}) to $V_{\a}(x)$ for $\a\le0$ and $x\in\interior{\Om}$, where $\interior{\Om}$ denotes the interior of $\Om$. Namely, if an integral $\displaystyle \int_\Om\omega$ blows up on a subset $X$ of $\Om$, we expand $\displaystyle \int_{\Om\setminus N_\e(X)}\omega$ in a series of $\frac1\e$, where $N_\e(X)$ $(\e>0)$ is an $\e$-tubular neighbourhood of $X$, and take the constant term (\cite{O3}). 
In the case of $V_\a(x)$ $(\a\le0, x\in\interior\Om)$, we have $X=\{x\}$ and the series has only one divergent term which depends only on $\a$ and $\n$. 
Thus we can obtain a $1$-parameter family of potentials, which we denote by $V_\Om^{(\a)}$ ($\a\in\RR$), with $V_\Om^{(\a)}=V_\a$ when $\a>0$ or $x\not\in\Om$. 

In particular, when $\Om$ is convex and $x\in\interior\Om$,  $V_\Om^{(\a)}(x)$ can be expressed as 
\begin{equation}\label{first_formula}
V_\Om^{(\a)}(x)=\frac1\a\int_{S^{\n-1}}\left(\rhosub{\Om_{-x}}(v)\right)^\a\,d\sigma(v) \hspace{1cm}(\a\ne0),
\end{equation}
where $\sigma$ is the standard Lebesgue measure of $S^{\n-1}$ and $\rhosub{\Om_{-x}}:S^{\n-1}\to\RR_{>0}$ is a {\em radial function} of $\Om_{-x}=\{y-x\,|\,y\in\Om\}$ given by $\mbox{\large $\rho  $}_{\Om_{-x}}(v)=\sup\{a\ge0\,|\,x+av\in \Om\}$. 
Thus $V_\Om^{(\a)}(x)$ coincides with the {\em dual mixed volume} $\tilde V_\alpha$ as introduced by Lutwak (\cite{Lu75,Lu88}) up to multiplication by a constant. 
%

We show basic properties and give boundary integral expressions using Stokes theorem, which plays an important role in computing the derivatives and hence Laplacians. 

Next we study points where the extreme values $Mm^{(\a)}(\Om)$ of $V_\Om^{(\a)}$ are attained. 
We will call such points the {\em $r^{\an}$-centers of $\Om$}. 
To be precise, it is a point that gives the minimum value of $V_\Om^{(\a)}$ 
when $\a>\n$, the maximum value of $V_\Om^{(\a)}$ when $0<\a<\n$, and the maximum value of $V_\Om^{(\a)}$ in $\interior\Om$ when $\a\le0$. 
When $\a=\n$, as $V_\Om^{(\n)}(x)$ is constantly equal to the volume of $\Om$, we define $r^0$-center by a point where the maximum value of the log potential is attained. 
We show that any compact set has an $r^{\an}$-center for any $\a$. 

For example, the centroid or the center of mass $x_G$ of $\Om$ is an $r^2$-center by the following reason. 
As $x_G$ is given by 
$\displaystyle x_G={\displaystyle \int_{\Om}y\,d\mu(y)}\Big/{\displaystyle \int_{\Om}1\,d\mu(y)}\,,$ 
or equivalently, by 
$\displaystyle \int_{\Om}(x_G-y)\,d\mu(y)=0\,,$
it can be characterized as the unique critical point of the map $V_\Om^{(\n+2)}:\RR^\n\ni x\mapsto\int_{\Om}{|x-y|}^2\,d\mu(y)\in\RR.$ 

Another example is the ``{\sl illuminating center}'' of a triangle, which was recently introduced by Katsuyuki Shibata (\cite{S}). 
It is a point that maximizes the total brightness of a triangular park obtained by a light source on that point. 
In our language, it is an $r^{-2}$-center. 
This is one of the motivations of our research. 

When $\Om$ is convex, an $r^\an$-center coincides with the {\em radial center of order $\a$}, which was introduced for $0<\a\le1$ in \cite{M1}. 

Furthermore, using the {\sl moving plane method} from analysis, we introduce a region whose complement has no chance to have any $r^{\an}$-center, which we call the {\em minimal unfolded region of $\Om$}. 

The $r^{\an}$-center is not necessarily unique. For example, a disconnected region has at least two $r^{\an}$-centers for small $\a$. 
We show that the $r^{\an}$-center is unique for any compact set if $\a\ge \n+1$, which includes the case of centroids, and for any {\sl convex} sets if $\a\le 1$. 
The latter is a consequence of a theorem in \cite{M1} when we use \eqref{first_formula}. 

We also show that among bodies with the same volume, only balls can give maximum (or minimum according to $\a$) of the extreme value $Mm^{(\a)}$ of $V_\Om^{(\a)}$ for any $\a$, and the same statement holds for the integration of the potential $V_\Om^{(\a)}$ on $\Om$ if $\a>0$, i.e. when the renormalization is not needed. 

Finally we study the asymptotic behavior of $r^\an$-centers as $\a$ goes to $\pm\infty$. 
We show that as $\a$ goes to $+\infty$ an $r^\an$-center approaches a point where the infimum of a map $\mathbb R^\n\ni x\mapsto \max_{y\in\Om}|y-x|\in\RR$ is attained, which we call a min-max point.   
On the other hand, if a point where the supremum of a map $\mathbb R^\n\ni x\mapsto \min_{y\in\overline{{\Om}^c}}|y-x|\in\RR$ is attained, which we call a max-min point, is uniquely determined, then the set of $r^\an$-centers converges to a singleton consisting of the max-min point as $\a$ goes to $-\infty$ with respect to the Hausdorff distance. 

\medskip
{\bf Assumption and notation}. 

We always assume that our $\Om$ is a body, that is, a compact set which is a closure of its interior, and that it has a piecewise $C^1$ boundary $\partial\Om$. 
We may require additional conditions, $(\ast\ast)$ in Lemma \ref{lemma_divergence_boundary} or $(\ast\ast\ast)$ in Lemma \ref{lemma_partial_differentiability}. 
Even with those, compact $\n$-dimensional submanifolds and polyhedra satisfy our conditions. 

We assume that the dimension $\n$ is greater than or equal to $2$ except in example \ref{example_pair_of_intervals}. 

\smallskip
Throughout the paper, $\interior X$ and $X^c$ denote the interior and complement of $X$. 
We denote the standard Lesbegue measure of $\RR^\n$ by $\mu$, and that of $\partial\Om$ and other $(\n-1)$-dimensional spaces like $S^{\n-1}$ by $\sigma$. 

\medskip
{\bf Acknowledgement}. 
The author would like to express his deep gratitude to Katsuyuki Shibata who informed him of the illuminating center of a triangle, and to his colleague Kazuhiro Kurata for many helpful suggestions in the case when $\n=2$ and $\a=0$, and for informing the author of the moving plane method. 
The author also thanks the referee deeply for careful reading and many invaluable suggestions. 

%
\section{Renormalized $r^{\an}$-potential}
\subsection{Definition}
Suppose $\a\le0$ and $x\in\interior\Om$. 
Fix a positive constant $R$ with $R<\textrm{dist}(x,\partial\Om)$. 
If we denote an $\n$-ball with center $x$ and radius $r>0$ by $B_r(x)$, we have 
\setlength\arraycolsep{1pt}
\[\begin{array}{rcl}
\displaystyle \int_{B_R(x)\setminus {B}_\e(x)}{|x-y|}^{\an}\,d\mu(y)
%
&=&\displaystyle \left\{
\begin{array}{ll}
\displaystyle A(S^{\n-1})\log\frac{R}{\e} & \hspace{0.5cm}\mbox{\rm if }\>\a=0\,,\\[1mm]
\displaystyle \frac{A(S^{\n-1})}{\a}\left(R^{\a}-\e^{\a}\right)& \hspace{0.5cm}\mbox{\rm if }\>\a<0\,,
\end{array}
\right.
\end{array}\]
where $A(S^{\n-1})$ is the volume of the $(n-1)$-dimensional unit sphere $S^{\n-1}$: 
\[
\displaystyle A(S^{\n-1})=\frac{2\pi^{\frac{\n}{2}}}{\varGamma(\frac{\n}2)}
=\left\{
\begin{array}{ll}
\displaystyle \frac{2\pi^{k+1}}{k!} & \hspace{0.5cm}\mbox{\rm if }\>\n-1=2k+1\,,\\
\displaystyle \frac{2^{k+1}\pi^{k}}{(2k-1)(2k-3)\cdots1} & \hspace{0.5cm}\mbox{\rm if }\>\n-1=2k\,.
\end{array}
\right.
\]
As $\Om\setminus {B}_\e(x)=(\Om\setminus B_R(x))\,\cup\, (B_R(x)\setminus {B}_\e(x))$ and $\Om\setminus B_R(x)$ is independent of $\e$, if we expand $\int_{\Om\setminus {B}_\e(x)}{|x-y|}^{\an}\,d\mu(y)$ in a series in $\frac1\e$, we have only one divergent term which does not depend on the point $x$ and the set $\Om$. 
Thus we are lead to the following defintion: 
\begin{definition}\label{def_renormalization_Riesz} \rm 
A compact set which is a closure of its interior is called a {\em body} (\cite{Ga}). 
We always assume that $\Om$ is a body with a piecewise $C^1$ boundary $\partial\Om$ in what follows. 

(1) Let $\Om$ be a body in $\RR^\n$ with a piecewise $C^1$ boundary. 
%
Define $V^{(\a)}_\Om$
by the following. 

\begin{enumerate}
\item[(i)] When $\a\le0$ and $x\in\interior\Om$, define $V^{(\a)}_\Om(x)$ by 
\begin{equation}\label{def_V_a}
\begin{array}{rcl}
V^{(\a)}_\Om(x)
&=&\displaystyle \lim_{\e\to+0}\left(\int_{\Om\setminus B_\e(x)}{|x-y|}^\an\,d\mu(y)-\frac{A(S^{\n-1})}{-\a}\cdot\frac1{\e^{-\a}}\right)
\end{array}
\end{equation}
for $\a<0$, and 
\begin{equation}\label{def_V_0}
\begin{array}{rcl}
V^{(0)}_\Om(x)
&=&\displaystyle \lim_{\e\to+0}\left(\int_{\Om\setminus B_\e(x)}{|x-y|}^{-\n}\,d\mu(y)-A(S^{\n-1})\log\frac1\e\right)
\end{array}
\end{equation}
for $\a=0$, 
and call them the {\em renormalization of $r^\an$-potential of $\Om$}. 

\item [(ii)] When $\a>0$ or $x$ is not in $\Om$ we put 
$$V^{(\a)}_\Om(x)=\int_{\Om}{|x-y|}^\an\,d\mu(y)\,.$$
%
\end{enumerate}

(2) Let $-\Om$ denote the same space $\Om$ with the reversed orientation. 
Define $V^{(\a)}_{-\Om}(x)=-V^{(\a)}_\Om(x)$. 

Let $\Om_j$ $(1\le j \le k)$ be a body with a piecewise $C^1$ boundary 
so that $\interior\Om_i\cap\interior\Om_j=\emptyset$ $(i\ne j)$. 
Then define $V^{(\a)}_{\cup_j\Om_j}(x)=\sum_jV^{(\a)}_{\Om_j}(x)$. 
\end{definition}

\begin{remark}\label{remark_log}\rm
In our study, if we put $\frac1ar^a=\log r$ $(r>0)$ when $a=0$ in our formulae for general cases, we obtain those for special cases, although $\log r=\lim_{a\to0}\frac1a(r^a-1)$ in fact. 
Such examples can be found in  
\eqref{def_V_a} and \eqref{def_V_0}, \eqref{f_Va_epsilon} and \eqref{f_V0_epsilon}, \eqref{f_V_ball}, \eqref{V_alpha_boundary} and \eqref{V0_boundary}, \eqref{V_alpha_\n=2_contour} and \eqref{V0_\n=2_contour}, \eqref{eq_div_V_alpha_boundary} and \eqref{eq_div_V_0_boundary}, \eqref{V^a_int_S^m-1} and \eqref{V^0_int_S^m-1}, 
and Theorem \ref{thm_radial_general}. 
\end{remark}

The renormalization in the case when $m=2$ and $\a=-2$ was given by Auckly and Sadun with more generality in \cite{AS}, where they studied surface energy. 

We remark that when $\a=\n$ we have $V^{(\n)}_{\Om}(x)=\textrm{Vol}(\Om)$ for any $x$. But when we are concerned with $r^{\an}$-centers, it is natural to use the log potential for $V^{(\n)}_{\Om}$. 

\subsection{Basic properties}

\begin{proposition}
\begin{enumerate}
\item The map $V_{\Om}^{(\a)}$ inherits symmetry from $\Om$, i.e. if $g$ is an isometry of $\RR^\n$ then $V_{g\cdot\Om}^{(\a)}(g\cdot x)=V_{\Om}^{(\a)}(x)$. 
\item Under a homothety $\RR^\n\ni x\mapsto k x\in \RR^\n$ $(k>0)$, 
\[
V_{k\Om}^{(\a)}(kx)=\left\{
\begin{array}{ll}
\displaystyle k^{\a}V_{\Om}^{(\a)}(x) & \hspace{0.3cm}\mbox{if $\a\ne0$},\\
\displaystyle V_{\Om}^{(0)}(x)+A(S^{\n-1})\log k& \hspace{0.3cm}\mbox{if $\a=0$ and $x\in\interior\Om$,}\\[1mm]
\displaystyle V_{\Om}^{(0)}(x)& \hspace{0.3cm}\mbox{if $\a=0$ and $x\not\in\Om$.}
\end{array}
\right.
\]
\end{enumerate}
\end{proposition}

\begin{figure}[htbp]
\begin{center}
\includegraphics[width=.35\linewidth]{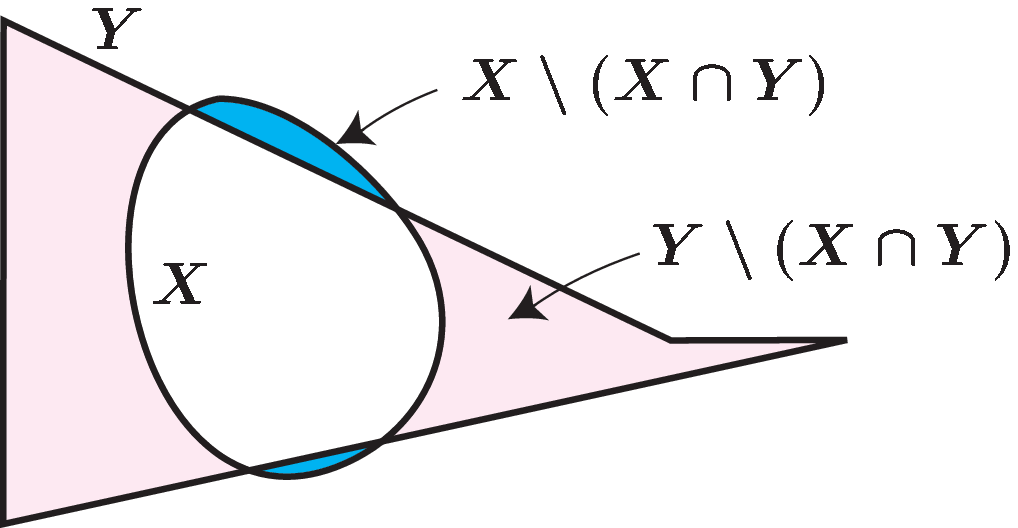}
\caption{$X-Y$. The dark (light) part with positive (resp. negative) orientation.}
\label{X-Y}
\end{center}
\end{figure}

\begin{proof} 
(2) follows directly from the definition, or from Theorem \ref{prop_boundary_integral_formula}. \end{proof}
\begin{lemma}\label{lemma_X-Y}
Let $X$ and $Y$ be subsets of $\RR^\n$. 
Define 
$$X-Y=(X\setminus(X\cap Y))\cup-(Y\setminus(X\cap Y)),$$
where the second term is equipped with the reverse orientation (figure \ref{X-Y}). 

If $x\in\interior{\Om}_1\cap\interior{\Om}_2$ or $x\in\Om_1^{\>c}\cap\Om_2^{\>c}$ then 
\[V^{(\a)}_{\Om_1-\Om_2}(x)=V^{(\a)}_{\Om_1}(x)-V^{(\a)}_{\Om_2}(x)\]
for any $\a$. 
\end{lemma}
\begin{proof}
This is because 
\[
\Om_1-\Om_2=(\Om_1\setminus B_\e(x))-(\Om_2\setminus B_\e(x))
\]
for sufficiently small $\e$ when $x\in\interior{\Om}_1\cap\interior{\Om}_2$. 
\end{proof}

We give elementary statements on the renormalization process of $V^{(\a)}_\Om$. 
\begin{proposition}
Suppose $x\in\interior\Om$. 
\begin{enumerate}
\item We have 
\begin{eqnarray}
V^{(\a)}_\Om(x)&=&\displaystyle \int_{\Om- B_\e(x)}{|x-y|}^\an\,d\mu(y)-\frac{A(S^{\n-1})}{-\a}\cdot\frac1{\e^{-\a}} \hspace{0.5cm}(\a\ne0),\label{f_Va_epsilon} \\
V^{(0)}_\Om(x)&=&\displaystyle \int_{\Om- B_\e(x)}{|x-y|}^{-\n}\,d\mu(y)-A(S^{\n-1})\log\frac1\e \hspace{0.86cm}(\a=0)\, \label{f_V0_epsilon}
\end{eqnarray}
for any $\e>0$. 
In particular, \eqref{def_V_a} and \eqref{def_V_0} hold even without taking the limit if $\e<\textrm{\rm dist}(x,\partial\Om)$. 
\item If $\a<0$ we have 
\begin{equation}
V^{(\a)}_\Om(x)= -\int_{\Om^{c}}{|x-y|}^{\an}\,d\mu(y) \hspace{1cm}(\a<0) \label{f_Va_complement}
\end{equation}
Thus the renormalization can be obtained by taking the complement of the domain when $\a<0$. 
\end{enumerate}
\end{proposition}

We remark that \eqref{f_V0_epsilon} is equivalent to 
\begin{equation}\label{f_V0_R}
V^{(0)}_\Om(x)=A(S^{\n-1})\log R-\int_{B_R(x)-\Om}{|x-y|}^{-\n}\,d\mu(y) 
\end{equation}
for any $R>0$. 
\begin{proof}
First note that the (renormalized) potential of a ball at the center is given by 
\begin{equation}\label{f_V_ball}
V^{(\a)}_{B_r(x)}(x)=\left\{
\begin{array}{ll}
\displaystyle \frac{A(S^{\n-1})}{\a}\cdot r^{\a} &\hspace{0.5cm} (\a\ne0),\\[3mm]
\displaystyle A(S^{\n-1})\log r &\hspace{0.5cm} (\a=0).
\end{array}\right.
\end{equation}

(1) By Lemma \ref{lemma_X-Y} we have 
$$
V^{(\a)}_\Om(x)=V^{(\a)}_{\Om-B_\e(x)}(x)+V^{(\a)}_{B_\e(x)}(x),
$$
which implies \eqref{f_Va_epsilon} and \eqref{f_V0_epsilon}. 

(2) Suppose $\a<0$ and $x\in\interior\Om$. Direct calculation shows 
$$
\int_{\RR^\n\setminus {B}_\e(x)}{|x-y|}^\an\,d\mu(y)
=-\frac{A(S^{\n-1})}{\a}\cdot{\e^{\a}}
=-V^{(\a)}_{B_\e(x)}(x),
$$
and therefore we have 
\[\begin{array}{rcl}
V^{(\a)}_\Om(x)&=&\displaystyle \int_{\Om- B_\e(x)}{|x-y|}^\an\,d\mu(y)-\int_{\RR^\n- {B}_\e(x)}{|x-y|}^\an\,d\mu(y) \\[4mm] 
&=&\displaystyle -\int_{\Om^{c}}{|x-y|}^{\an}\,d\mu(y)\,.
\end{array}\]
\end{proof}
\begin{corollary}\label{cor_Omega_1_subset_Omega_2}
Suppose $\Om_1\subset\Om_2$. 
If $\a>0$ or $x\in\interior{\Om_1}$ or $x\in\Om_2^c$ then $V_{\Om_1}^{(\a)}(x)\le V_{\Om_2}^{(\a)}(x)$. 
\end{corollary}

\begin{lemma}\label{lem_boundary} 
Let $V_{\Om}^{(\a)}\big|_{\interior\Om}$ and $V_{\Om}^{(\a)}\big|_{\Om^{c}}$ denote the restrictions of $V_{\Om}^{(\a)}$ to the interior and the complement of $\Om$ respectively. 
\begin{enumerate}
\item If $\a>0$ then $V_{\Om}^{(\a)}>0$ on $\RR^\n$. 
\item $V_{\Om}^{(\a)}\big|_{\Om^{c}}>0$ for any $\a$. 
If $\a<\n$ then $\displaystyle \lim_{|x|\to+\infty}V_{\Om}^{(\a)}(x)=0$. 
\item If $\a<0$ then $V_{\Om}^{(\a)}\big|_{\interior\Om}<0$. 
\end{enumerate}
\end{lemma}
\begin{proof}
(3) follows from the formula \eqref{f_Va_complement}. 
\end{proof}
%
\subsection{Boundary integral expressions and derivatives}
%
\begin{theorem}\label{prop_boundary_integral_formula}
If $x\not\in\partial\Om$ 
then $V_{\Om}^{(\a)}(x)$ can be expressed by the boundary integral on $\partial\Om$ as 
\begin{eqnarray}
V_{\Om}^{(\a)}(x)&=&\displaystyle \displaystyle \frac1{\a}\int_{\partial\Om} {|x-y|}^{\an}(y-x)\cdot n\,d\sigma (y) \hspace{0.5cm} (\a\ne0),\label{V_alpha_boundary}\\[1mm]
V_{\Om}^{(0)}(x)&=&\displaystyle \displaystyle \int_{\partial\Om} \frac{\log|x-y|}{{|x-y|}^\n}\,(y-x)\cdot n\,d\sigma (y) \hspace{0.9cm} (\a=0), \label{V0_boundary}
\end{eqnarray}
where $n$ is a unit outer normal vector to $\partial\Om$ at $y$ and $\sigma$ denotes the standard Lebesgue measure of $\partial\Om$. 
\end{theorem}

Especially, when $\n=2$, $V_{\Om}^{(\a)}(x)$ $(\a\ne0)$ can be expressed by the contour integral along $\partial\Om$ by 
\begin{equation}\label{V_alpha_\n=2_contour}
\begin{array}{rcl}
V_{\Om}^{(\a)}(x)&=&\displaystyle \frac1{\a^2}\oint_{\partial\Om} \left(\nabla_y{|x-y|}^\a\right)\cdot n\,ds \hspace{2.0cm} (\a\ne0)\\[4mm]
&=&\displaystyle \frac1{\a}\oint_{\partial\Om} {|x-y|}^{\a-2}\left((y_1-x_1)\,dy_2-(y_2-x_2)\,dy_1\right)\,  
\end{array}
\end{equation}
where $s$ is the arc-length of $\partial \Om$, and $V_{\Om}^{(0)}(x)$ by 
\begin{equation}\label{V0_\n=2_contour}
\begin{array}{rcl}
V_{\Om}^{(0)}(x)&=&\displaystyle \oint_{\partial\Om} \log|x-y|\,(\nabla_y\log|x-y|)\cdot n\,ds \\[4mm]
&=&\displaystyle \oint_{\partial\Om} \frac{\log|x-y|}{{|x-y|}^2}\left((y_1-x_1)\,dy_2-(y_2-x_2)\,dy_1\right)\,.
\end{array}
\end{equation}

The author thanks Kazuhiro Kurata for informing him of the first formula of \eqref{V0_\n=2_contour}. 

We remark that the $1$-form in \eqref{V0_\n=2_contour} can be expressed as $-\log|x-y|\,\Im\mathfrak{m}\frac{dy}{x-y}$ if we consider $x$ and $y$ as points in $\mathbb{C}$. 
\begin{proof} Put $r=|x-y|$. 

(1-i) Suppose $\a<0$ and $x\in\interior\Om$. 
Note that 
\begin{equation}\label{eq_div_V_alpha_boundary}
\textrm{div}_y\left(\frac1{\a}r^{\an}(y-x)\right)
=r^\an.
\end{equation}
Therefore, for enough small $\e>0$, 
\[\begin{array}{rcl}
\displaystyle \int_{\Om\setminus B_\e(x)}r^\an\,d\mu(y)
&=&\displaystyle \frac1{\a}\int_{\partial\Om} r^{\an}(y-x)\cdot n\,d\sigma (y)
-\frac1{\a}\int_{\partial B_\e(x)} r^{\an}(y-x)\cdot n\,d\sigma (y) \\[4mm]
&=&\displaystyle \frac1{\a}\int_{\partial\Om} r^{\an}(y-x)\cdot n\,d\sigma (y)
-\frac{A(S^{\n-1})\,\e^\a}{\a}\,.
\end{array}\]
Then the formula \eqref{f_Va_epsilon} implies \eqref{V_alpha_boundary}. 

\smallskip
(1-ii) The case when $x\not\in\Om$ or when $\a>0$ 
can be proved in the same way if we forget the renormalization term. 

\smallskip
(2) The case when $\a=0$ can be proved in the same way using 
\begin{equation}\label{eq_div_V_0_boundary}
\textrm{div}_y\left((r^{-\n}{\log r})(y-x)\right)
=r^{-\n}
\end{equation}
and \eqref{f_V0_epsilon}
\end{proof}

\begin{proposition}\label{prop_derivative_boundary} 
{\rm (1)} 
If $x=(x_1, \ldots, x_\n)\not\in\partial\Om$ 
then $\displaystyle \frac{\partial V_{\Om}^{(\a)}}{\partial x_j}(x)$ can be expressed by the boundary integral on $\partial\Om$ as 
\begin{equation}\label{formula_derivative_boundary}
\frac{\partial V_{\Om}^{(\a)}}{\partial x_j}(x)=-\int_{\partial\Om}{|x-y|}^{\an}\,e_j\cdot n\,d\sigma (y)
\end{equation}
for any $j$ $(1\le j\le \n)$ and for any $\a$, 
where $n$ is a unit outer normal vector to $\partial\Om$ at $y$, $e_j$ is the $j$-th unit vector of $\RR^\n$, and $\sigma$ denotes the standard Lebesgue measure of $\partial\Om$. 

\smallskip
{\rm (2)} If $\a>1$ then \eqref{formula_derivative_boundary} holds and becomes continuous on entire $\RR^\n$. 
Therefore $V^{(\a)}_\Om$ is of class $C^1$ on $\RR^\n$ if $\a>1$.  
\end{proposition}
\begin{proof}
(1) Put $r=|x-y|$. 
Note that 
$$
\frac{\partial \, r^{\an}}{\partial x_j}=-\frac{\partial \, r^{\an}}{\partial y_j}=-\textrm{div}_y\left(r^\an \,e_j\right)
$$

(i) Assume $x\not\in\Om$. 
Then 
\[
\begin{array}{rcl}
\displaystyle \frac{\partial V_{\Om}^{(\a)}}{\partial x_j}(x)
&=&\displaystyle \int_\Om\frac{\partial \, r^{\an}}{\partial x_j}\,d\mu(y)
=-\int_\Om\textrm{div}_y\left(r^\an \,e_j\right)\,d\mu(y)
=\displaystyle -\int_{\partial\Om}r^\an \,e_j\cdot n\,d\sigma (y) 
\end{array}
\]
as is required. 

\smallskip
(ii) Assume $x\in\interior\Om$. 
Put $x_0=x$ and fix it. 
Take a positive number $\e_0$ with $\e_0<\textrm{dist}(x_0,\partial\Om)$. 
Suppose $x'\in B_{\frac{\e_0}2}(x_0)$. 
Then 
\[V_{\Om}^{(\a)}(x')-V_{B_{\e_0}(x_0)}^{(\a)}(x')=\int_{\Om\setminus B_{\e_0}(x_0)}|x'-y|^\an\,d\mu(y)\,.\]
Therefore, 
\[\begin{array}{rcl}
\displaystyle \frac{\partial V_{\Om}^{(\a)}}{\partial x_j}(x')-\frac{\partial V_{B_{\e_0}(x_0)}^{(\a)}}{\partial x_j}(x')
&=&\displaystyle \int_{\Om\setminus B_{\e_0}(x_0)} \frac{\partial }{\partial x_j}|x'-y|^\an\,d\mu(y)\\[4mm]
&=&\displaystyle -\int_{\Om\setminus B_{\e_0}(x_0)} \textrm{div}_y\left(|x'-y|^\an \,e_j\right)d\mu(y)\\[4mm]
&=&\displaystyle -\int_{\partial\Om}|x'-y|^\an \,e_j\cdot n\,d\sigma(y)
+\int_{\partial B_{\e_0}(x_0)}|x'-y|^\an \,e_j\cdot n\,d\sigma(y)\,,
\end{array}\]
where $\sigma$ is the standard Lesbegue measure of $\partial B_{\e_0}(x_0)$. 

If we put $x'=x_0$, since 
\[
\frac{\partial V_{B_{\e_0}(x_0)}^{(\a)}}{\partial x_j}(x_0)=0\>\>\mbox{ and }\>\>
\int_{\partial B_{\e_0}(x_0)}|x_0-y|^\an \,e_j\cdot n\,d\sigma(y)=0
\]
by the symmetry, we obtain \eqref{formula_derivative_boundary}. 

\medskip
(2) Suppose $\a>1$. 
Since $\a-\n>-(\n-1)=-\dim\partial\Om$, the right hand side of \eqref{formula_derivative_boundary} is well-defined even if the point $x$ is on the boundary $\partial\Om$, and is continuous on $\RR^\n$. 
\end{proof}
\begin{remark}\label{geometric_proof}\rm 
There are two more proofs of the above proposition. 

(i) Put $\Om_{-he_j}=\{y-he_j\,|\,y\in\Om\}$. Then (see Figure \ref{Omega_-he1}) 
$$
\begin{array}{rcl}
\displaystyle \frac{\partial V_{\Om}^{(\a)}}{\partial x_j}(x)
&=&\displaystyle \lim_{h\to0}\frac{V_{\Om}^{(\a)}(x+he_j)-V_{\Om}^{(\a)}(x)}h \\[4mm]
&=&\displaystyle \lim_{h\to0}\frac{V_{\Om_{-he_j}}^{(\a)}(x)-V_{\Om}^{(\a)}(x)}h \\[4mm]
&=&\displaystyle -\lim_{h\to0}\frac1h \int_{\Om-\Om_{-he_j}}|x-y|^\an\,d\mu(y)\\[4mm]
&=&\displaystyle -\int_{\partial\Om}{|x-y|}^{\an}\,e_j\cdot n\,d\sigma (y).
\end{array}
$$

\begin{figure}[htbp]
\begin{center}
\includegraphics[width=.6\linewidth]{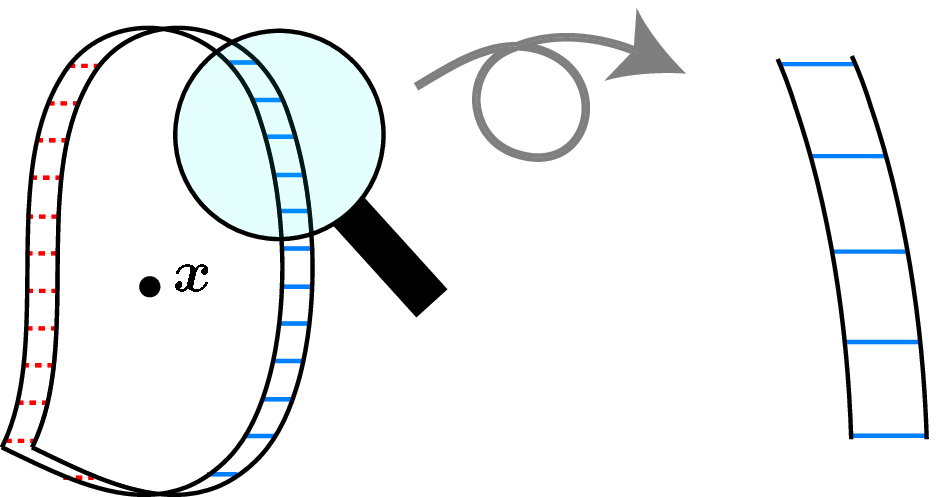}
\caption{$\Om-\Om_{-he_1}$ $(h>0)$. The right (left) strip with positive (resp. negative) orientation. }
\label{Omega_-he1}
\end{center}
\end{figure}

\medskip
(ii) The proposition can also be proved by the boundary integral formulae in Theorem \ref{prop_boundary_integral_formula}, Stokes' theorem, and Lemma \ref{lemma_X-Y}. 
\end{remark}

\begin{corollary}\label{cor_second_partial_derivatives}
If $\a>2$ 
then the second derivative $\frac{\partial^2 V_{\Om}^{(\a)}}{\partial x_j{}^2}$ is given by 
\begin{eqnarray}
\displaystyle \frac{\partial^2 V_{\Om}^{(\a)}}{\partial x_j{}^2}(x)
&=&\displaystyle -(\a-\n)\int_{\partial\Om}{|x-y|}^{\an-2}(x_j-y_j)\,e_j\cdot n\,d\sigma (y)  \label{f_second_partial_derivative_boundary} \\
&=&\displaystyle (\a-\n)\int_{\Om}{|x-y|}^{\an-4}\left((\an-2)(x_j-y_j)^2+{|x-y|}^2\right)d\mu (y) \label{f_second_partial_derivative_Omega_bis}\\
&=&\displaystyle (\a-\n)\int_{\Om}{|x-y|}^{\an-4}\left((\an-1)(x_j-y_j)^2+\sum_{i\ne j}(x_i-y_i)^2\right)\!d\mu (y),\hspace{0.3cm}{\phantom{a}}
\label{f_second_partial_derivative_Omega}
\end{eqnarray}
and is continuous on $\RR^\n$. 

Furthermore, for any $\a$,  
\eqref{f_second_partial_derivative_boundary} holds on $(\partial\Om)^c$ and \eqref{f_second_partial_derivative_Omega_bis} and \eqref{f_second_partial_derivative_Omega} hold on $\Om^{\,c}$. 

\end{corollary}
\begin{proof}
Put $r=|x-y|$. 

The absolute value of the integrand of \eqref{f_second_partial_derivative_boundary} is bounded above by $r^{\an-2}\cdot r=r^{\an-1}$. 
If $\a>2$ then $\an-1>-m+1=-\dim(\partial\Om)$, and therefore the right hand side of \eqref{f_second_partial_derivative_boundary} is well-defined even if a point $x$ belongs to $\partial\Om$. 

The absolute value of the integrand of \eqref{f_second_partial_derivative_Omega} is bounded above by $(|\an-1|+1)\,r^{\an-2}$. 
If $\a>2$ then $\an-2>-\n$, and therefore \eqref{f_second_partial_derivative_Omega} is well-defined. 

The equality of the right hand side of \eqref{f_second_partial_derivative_boundary} and \eqref{f_second_partial_derivative_Omega_bis} is obvious.
\end{proof}
%
\subsection{Continuity and asymptotic behavior}
%
The same geometric argument as in Remark \ref{geometric_proof} (i) implies the continuity of $V_{\Om}^{(\a)}$. 
\begin{proposition}\label{prop_continuity_homothety} 
When $\a>0$ the map $V_{\Om}^{(\a)}:\RR^\n\to\RR$ is continuous.  
When $\a\le0$, the restrictions of $V_{\Om}^{(\a)}$ to $\interior{\Om}$ and $\Om^{\,c}$ are continuous. 
\end{proposition}

We remark that when $0<\a<\n$ the continuity is nothing but that of the Riesz potential, 
and that when $\a\le0$ it follows from the boundary integral expression in Theorem  \ref{prop_boundary_integral_formula}. 
%
\begin{lemma}\label{lemma_divergence_boundary}
Suppose a point $z$ on the boundary $\partial\Om$ satisfies Poincar\'e's condition outside, namely, 

\smallskip
$(\ast\ast)$ There are positive constants $\e$ and $\theta$ such that 
$B_{\e}(z)\cap\Om^{\,c}$ contains a cone of revolution with vertex $z$ and cone angle $\theta$. 

\begin{enumerate}
\item If $\a\le0$ then $V_{\Om}^{(\a)}\big|_{\interior\Om}(x)$ goes to $-\infty$ uniformly as $x\in\interior\Om$ approaches $z$. 
\item Assume that $\partial\Om$ satisfies Poincar\'e's condition inside. 
If $\a\le0$ then $V_{\Om}^{(\a)}\big|_{\Om^{c}}(x)$ goes to $+\infty$ uniformly as $x\in\Om^{\,c}$ approaches $z$. 
\end{enumerate}
\end{lemma}
\begin{proof}
We only prove (1) in what follows as (2) can be proved in the same way. 

\smallskip
We may assume $\e\le\frac12$. 
Let $C$ be the volume of a subset of the unit sphere $S^{\n-1}$ which is inside a cone of revolution with center the origin and cone angle $\theta$. 
Suppose $x\in\interior\Om$ satisfies $|x-z|<\delta$, where $\delta\le\frac12$. 

(i) Assume $\a<0$. 
Since $\delta+r\le1$ if $0<r<\frac12$, the formula \eqref{f_Va_complement} implies 
\[\begin{array}{rcl}
V_{\Om}^{(\a)}(x)&\le&\displaystyle -\int_0^{\e}(\delta+r)^\an\cdot Cr^{\n-1}\,dr 
%
\le\displaystyle -C\int_0^{\e}{(\delta+r)}^{-\n}\cdot r^{\n-1}\,dr
=-C\int_1^{1+\frac{\e}\delta} \frac{(t-1)^{\n-1}}{t^\n}\,dt,
\end{array}\]
which implies $\displaystyle \lim_{|x-z|\to+0}V_{\Om}^{(\a)}(x)=-\infty$ as $\displaystyle \lim_{R\to+\infty}\int_1^R \frac{(t-1)^{\n-1}}{t^\n}\,dt=+\infty$. 

\smallskip
(ii) Assume $\a=0$. 
Put $R=\max\{1,\diam\}$, where $\diam$ is the diameter of $\Om$. 
Then $B_R(x)$ contains both $\Om$ and the cone in $\Om^{\,c}$ of the condition $(\ast\ast)$. 
Therefore, the formula \eqref{f_V0_R} implies $\displaystyle \lim_{|x-z|\to+0}V_{\Om}^{(0)}(x)=-\infty$ by a similar argument as above. 

\smallskip
The uniformness follows from the fact that the right hand sides of the above estimates do not depend on $x$. 
\end{proof}
\begin{corollary}
If $\a\le0$ then $V_{\Om}^{(\a)}$ cannot be extended to a continuous function on $\RR^\n$. 
\end{corollary}
\begin{lemma}\label{lemma_partial_differentiability}
Suppose $\a\le1$. 
Assume $\Om$ satisfies the following condition: 

\smallskip
$(\ast\ast\ast)$ For any boundary point $z$ there is a direction $v\in S^{\n-1}$ such that $\Om$ around that point can be obtained as a space under a graph in direction $v$ of a piecewise $C^1$ function on an $(\n-1)$ dimensional hyperplane. 
\smallskip

Let $v$ be a boundary point and $v$ a unit outer vector as above. 
Then the directional derivative of $V^{(\a)}_\Om$ in $v$ satisfies 
\[\lim_{h\to0}\frac{\partial \,V^{(\a)}_\Om}{\partial v}(z+hv)=-\infty\,.\]

\end{lemma}
First observe that the formula \eqref{formula_derivative_boundary} implies that the lemma above holds when $\partial\Om$ contains $B_R^{\n-1}(0)\times\{0\}$ for some $R>0$ with outer normal vector $e_\n$ and the point $z$ approaches $\partial\Om$ as $z=he_\n$ $(h\to0)$. The proof of a general case is just a technicality. 

\begin{proof} 
We may assume, after a motion of $\RR^\n$, that $z=0$ and $v=e_m$, namely, there are positive constants $\rho$ and $b$ and a piecewise $C^1$ function $f$ such that 
\[\begin{array}{l}
\left(B_{\rho}^{\n-1}(0)\times(-b,b)\right)\cap\Om\\[1mm]
=\{(x_1,\ldots,x_\n)\,|\,x_\n\le f(x_1, \ldots, x_{\n-1}), \>(x_1, \ldots, x_{\n-1})\in B_{\rho}^{\n-1}(0)\}\,.
\end{array}\]

Let $C'>0$ be the minimum of $e_\n\cdot n$ on $\left(B_{\rho}^{\n-1}(0)\times[-b,b]\right)\cap\partial\Om$. 
Let $\pi:\RR^\n\to\RR^{\n-1}$ be the projection in the direction of $e_\n$. 
There is a positive constant $C$ such that for any point $y$ in $\left(B_{\rho}^{\n-1}(0)\times(-b,b)\right)\cap\partial\Om$ there holds $|y|<C|\pi(y)|$. 
For any $\e>0$ if $0<|h|<\e$ then for any $y\in\partial\Om$ with $\pi(y)\in B_{\rho}^{\n-1}(0)\setminus B_\e^{\n-1}(0)$ 
we have 
$\displaystyle \frac{|he_m-y|}{|y|}\le\frac{|h|+|y|}{|y|}<2.$
Then if $0<|h|<\e$ then 
\[\begin{array}{l}
\displaystyle \int_{\left(B_{\rho}^{\n-1}(0)\times(-b,b)\right)\cap\partial\Om}
|he_\n-y|^{\an}\,e_\n\cdot n\,d\sigma (y)\\[4mm]
\displaystyle >\int_{\left((B_{\rho}^{\n-1}(0)\setminus B_\e^{\n-1}(0))\times(-b,b)\right)\cap\partial\Om}
|he_\n-y|^{\an}\,e_\n\cdot n\,d\sigma (y)\\[4mm]
\displaystyle >2^\an\int_{\left((B_{\rho}^{\n-1}(0)\setminus B_\e^{\n-1}(0))\times(-b,b)\right)\cap\partial\Om}
|y|^{\an}\,e_\n\cdot n\,d\sigma (y)\\[4mm]
\displaystyle \ge2^\an C^\an C'\int_{B_{\rho}^{\n-1}(0)\setminus B_\e^{\n-1}(0)}|y'|^{\an}\,d\mu_{\n-1} (y')\\[4mm]
\displaystyle =2^\an C^\an C' A(S^{\n-2})\int_\e^\rho r^\an\cdot r^{\n-2}\,dr\,,
\end{array}\]
where $\mu_{\n-1}$ is the standard Lebesgue measure of $\RR^{\n-1}$. 
As $\a-2\le-1$ the right hand side above blows up to $+\infty$ as $\e$ goes down to $+0$, and therefore, 
\[\begin{array}{rcl}
\displaystyle \frac{\partial V_{\Om}^{(\a)}}{\partial x_\n}(he_\n)&=&\displaystyle -\int_{\partial\Om}|he_\n-y|^{\an}\,e_\n\cdot n\,d\sigma (y)
%
\end{array}\]
goes to $-\infty$ as $h$ approaches $0$. 
\end{proof}

\medskip
Let us consider the continuity of $V_{\Om}^{(\a)}(x)$ with respect to $\a$. 
\begin{proposition}\label{prop_continuity} 
Fix a compact set $\Om$ and a point $x\in\RR^\n$. 
\begin{enumerate}
\item The map $\RR\ni\alpha\mapsto V_{\Om}^{(\a)}(x)\in\RR$ is continuous on $\RR_+$ for any $x$. 
It is further continuous on $\RR\setminus\{0\}$ if $x\in\interior\Om$, and on $\RR$ if $x\not\in\Om$. 
\item If $x\in\interior\Om$ then $\displaystyle \lim_{\a\to\pm0}V_{\Om}^{(\a)}(x)=\pm\infty$. 
\end{enumerate}
\end{proposition}
\begin{proof}
Put $\varPhi_{\Om}(\a)=V_{\Om}^{(\a)}(x)$. 

(1) Suppose $\a>0$. 
Let $\e$ be a positive number. 
Then \eqref{f_V_ball} implies that there is $r>0$ such that $\varPhi_{\Om\cap B_r(x)}(\a')<\e$ 
for any $\a'$ with $\frac\a2\le\a'$. 
On the other hand, as $\Om\setminus {(\Om\cap B_r(x))}^\circ$ is compact and $\textrm{dist}(x,\Om\setminus {(\Om\cap B_r(x))}^\circ)\ge r$, where ${(\Om\cap B_r(x))}^\circ$ denotes the interior of $\Om\cap B_r(x)$, there is $\delta>0$ ($\delta\le\frac\alpha2$) such that if $|\a'-\a|<\delta$ then 
$\big|\varPhi_{\Om\setminus {(\Om\cap B_r(x))}^\circ}(\a')-\varPhi_{\Om\setminus {(\Om\cap B_r(x))}^\circ}(\a)\big|<\e$ 
and hence 
\[\begin{array}{rcl}
\displaystyle \left|\varPhi_{\Om}(\a')-\varPhi_{\Om}(\a)\right|
&\le&\displaystyle \left|\varPhi_{\Om\setminus {(\Om\cap B_r(x))}^\circ}(\a')-\varPhi_{\Om\setminus {(\Om\cap B_r(x))}^\circ}(\a)\right| + \varPhi_{B_r(x)}(\a')+\varPhi_{B_r(x)}(\a) \le3\e,
\end{array}\]
which implies the first statement. 

\smallskip
The second statement follows from \eqref{V_alpha_boundary}, and the last one is obvious. 

\smallskip
(2) Suppose $x\in\interior\Om$. 
Let $r$ be a positive number with $r<\textrm{dist}(x,\partial\Om)$ and $\diam$ be the diameter of $\Om$. 
Then $B_r(x)\subset\Om\subset B_\diam(x)$. 

By Corollary \ref{cor_Omega_1_subset_Omega_2} and the formula \eqref{f_V_ball}, if $\a>0$ then 
\[V^{(\a)}_\Om(x)\ge V^{(\a)}_{B_r(x)}(x)=\frac{A(S^{\n-1})\,r^{\a}}{\a},\]
which implies $\displaystyle \lim_{\a\to+0}\varPhi_{\Om}(\a)=+\infty$, and if $\a<0$ then 
\[V^{(\a)}_\Om(x)\le V^{(\a)}_{B_\diam(x)}(x)=\frac{A(S^{\n-1})\,\diam^{\phantom{,\!}\a}}{\a},
\]
which implies $\displaystyle \lim_{\a\to-0}\varPhi_{\Om}(\a)=-\infty$.
\end{proof}

The asymptotic behavior of $V^{(\a)}_\Om$ as $\a$ goes to $\pm\infty$ will be studied in Lemma \ref{l_asympt_V}. 

\subsection{Laplacian}
%
\begin{theorem}
Let $\Om$ be a body in $\RR^\n$ with a piecewise $C^1$ boundary. 
If $\a>2$ or $x\not\in\partial\Om$ and $\a\ne2$ then 
\begin{equation}\label{laplacian_general}
\Delta V^{(\a)}_\Om(x)=(\a-2)(\an)V^{(\a-2)}_\Om(x).
\end{equation}
%
\end{theorem}
\begin{proof} 
First note that 
\begin{eqnarray}\label{laplacian_r}
\Delta_x{|x-y|}^\an=(\a-2)(\a-\n){|x-y|}^{\an-2}.
\end{eqnarray}

If $\a>2$ then \eqref{laplacian_general} follows from \eqref{f_second_partial_derivative_Omega_bis} or from \eqref{laplacian_r}. 

Suppose $x\not\in\partial\Om$. 
By formula \eqref{f_second_partial_derivative_boundary}, we have 
\begin{eqnarray}
\Delta V^{(\a)}_\Om(x)&=&\displaystyle -(\a-\n)\int_{\partial\Om}{|x-y|}^{\an-2}(x-y)\cdot n\,d\sigma (y).  \label{f_laplacian_1}
\end{eqnarray}

On the other hand, if $\a-2\ne0$ then by Theorem \ref{prop_boundary_integral_formula} we have 
\[V_{\Om}^{(\a-2)}(x)=\frac1{\a-2}\int_{\partial\Om} {|x-y|}^{\an-2}(y-x)\cdot n\,d\sigma (y). \]
%
\end{proof}

We remark that when $\a=2$ and $\n\ge3$, i.e. when $V^{(\a)}_\Om$ is (a constant times) the Newton potential, we have $\Delta V^{(2)}_\Om(x)=0$ if $x\in\Om^{\,c}$ and $\Delta V^{(2)}_\Om(x)=-(\n-2)A(S^{\n-1})$ if $x\in\interior\Om$. 

\begin{corollary}
Let $\Om$ be a body in $\RR^\n$ with a piecewise $C^1$ boundary. 
\begin{enumerate}
\item If $\a>\n$ then $\Delta V^{(\a)}_\Om>0$ on $\RR^\n$. 
\item If $\a=\n$ then $\Delta V^{(\n)}_\Om=0$ on $\RR^\n$. 
\item If $2<\a<\n$ $(\n\ge3)$ then $\Delta V^{(\a)}_\Om<0$ on $\RR^\n$. 
\item If $\n\ge3$ and $\a=2$ then 
$\Delta V^{(2)}_\Om<0$ on $\interior\Om$. 
\item If $\a<2$ then $\Delta V^{(\a)}_\Om>0$ on $\Om^{\,c}$ and $\Delta V^{(\a)}_\Om<0$ on $\interior\Om$. 
\end{enumerate}
\end{corollary}

The author thanks Kazuhiro Kurata for informing him that $\Delta V^{(0)}_\Om<0$ on $\interior\Om$ when $\n=2$. 
\begin{proof}
The formula \eqref{laplacian_general} and Lemma \ref{lem_boundary} imply (1), (3), and (5). 
The statement (2) follows from $V^{(\n)}_\Om(x)=\textrm{Vol}(\Om)$ $(\forall x)$. 
The statement (4) follows from the remark above. 
\end{proof}

\subsection{Convex geometric formulae and duality by an inversion in a sphere}
%
Let us recall some terminologies from convex geometry. 

A set $M$ in $\RR^\n$ is called {\em star-shaped at $0$} if for any point $p$ in $M\setminus\{0\}$ the line segment $\overline{0p}$ is contained in $M$. 

Suppose $M$ is star-shaped at $0$. 
The {\em radial function} of $M$ is a map $\rhosub{M}$ from $S^{\n-1}$ to $\RR_{\ge0}$ defined by $\rhosub{M}(v)=\sup\{a\ge0\,|\,av\in M\}$ 
($\rhosub{M}(v)$ is also denoted by $\mbox{\large $\rho $}(M,v)$).
When $\interior{M}\ni0$ the {\em dual mixed volume of $M$ of order $\alpha$ $(\alpha\in\RR)$} is given by 
$$
\tilde{V}_\alpha(M)=\int_{S^{\n-1}}{(\rhosub{M}(v))}^\alpha d\sigma(v).
$$
When $\alpha=\n$, $\tilde{V}_\n(M)=\n\textrm{Vol}(M)$. 
It was originally introduced for convex bodies in \cite{Lu75,Lu88} (see also \cite{Ga}, \cite{Sch}). 

We show that $V^{(\a)}_\Om(x)$ $(\alpha\ne0, \, x\in\interior\Om)$ coincides with the dual mixed volume of $\Om_{-x}$ of order $\alpha$ up to multiplication by a constant. 
\begin{lemma}\label{prop_int_geom}
Let $\Om$ be a body in $\RR^\n$ with a piecewise $C^1$ boundary. 
 %
Suppose $x\in\interior{\Om}$ and $\Om$ is star-shaped at $x$. 
Let $\rhosub{\Om_{-x}}$ 
be a radial function of $\Om_{-x}=\{y-x\,|\,y\in\Om\}$. 
Then 
\begin{eqnarray}
V^{(\a)}_\Om(x) &=&\displaystyle \frac1\a\int_{S^{\n-1}}\left(\rhosub{\Om_{-x}}(v)\right)^\a\,d\sigma(v) \hspace{1cm}(\a\ne0), \label{V^a_int_S^m-1} \\
V^{(0)}_\Om(x) &=&\displaystyle \int_{S^{\n-1}}\log\left(\rhosub{\Om_{-x}}(v)\right)\,d\sigma(v) \hspace{1cm}(\a=0). \label{V^0_int_S^m-1} 
\end{eqnarray}
\end{lemma}
\begin{proof}
Note that $\Om$ can be expressed as 
$$
\Om=\left\{x+r(v)\,v\,|\,0\le r(v)\le \rhosub{\Om_{-x}}(v)\,, \,v\in S^{\n-1}\right\}.
$$

Remark that 
$r^{\an}d\mu=r^{\a-1}dr\,d\sigma(v).$ 

If $\a<0$ then \eqref{f_Va_complement} implies 
$$
V^{(\a)}_\Om(x)
=-\int_{S^{\n-1}}\int_{\rho _{\Om_{-x}}(v)}^\infty r^{\a-1} dr\,d\sigma(v)
=\frac1\a\int_{S^{\n-1}}\left(\mbox{\large $\rho$}_{\Om_{-x}}(v)\right)^\a d\sigma(v).
$$

If $\a=0$ then \eqref{f_V0_epsilon} implies that if 
we use a sufficiently small $\e>0$ 
we have  
$$
V^{(\a)}_\Om(x)
=\displaystyle \int_{S^{\n-1}}\int_\e^{\rho_{\Om_{-x}}(v)} \frac{dr}{r}\,d\sigma(v)-A(S^{\n-1})\log\frac1\e
=\int_{S^{\n-1}}\log\left(\rhosub{\Om_{-x}}(v)\right)\,d\sigma(v).
$$

The case when $\a>0$ is obvious. 
\end{proof}
\begin{theorem}\label{thm_radial_general}
Let $\Om$ be a body in $\RR^\n$ with a piecewise $C^1$ boundary. 
Let $\ell_{x,v}$ $(x\in\RR^\n, v\in S^{\n-1})$ denote a half line $\{x+tv\,|\,t>0\}$, and $p_i$ be a point where $\ell_{x,v}$ intersects $\partial\Om$ transversally, which we denote by $p_i\in\ell_{x,v}\ti\partial\Om$. 
Define the signature at $p_i$ by $\mbox{\rm sgn}(p_i)=+1$ if $\ell_{x,v}$ cuts across $\partial\Om$ from the inside to the outside of $\Om$ as $t$ increases and $\mbox{\rm sgn}(p_i)=-1$ otherwise. Then if $\a>0$ or $x\not\in\partial\Om$ we have 
\begin{eqnarray}
V^{(\a)}_\Om(x)&=&\displaystyle \frac1\a\int_{S^{\n-1}}
\sum_{p_i\in\ell_{x,v}\ti\partial\Om}\mbox{\rm sgn}(p_i)\,{|p_i-x|}^\a
\,d\sigma(v) \hspace{0.5cm}(\a\ne0), \nonumber\\
V^{(0)}_\Om(x)&=&\displaystyle \int_{S^{\n-1}}
\sum_{p_i\in\ell_{x,v}\ti\partial\Om}\mbox{\rm sgn}(p_i)\log|p_i-x|
\,d\sigma(v) \hspace{0.5cm}(\a=0). \nonumber
\end{eqnarray}
\end{theorem}

This can be proved by cutting $\Om$ into several pieces so that the intersection of every piece and any half line starting from $x$ is connected (possibly empty). 
\begin{corollary}\label{thm_inversion}
Suppose $x$ is a point with $x\not\in\partial\Om$. 
Let $I_x$ be an inversion in a unit sphere with center $x$. 
Let $\Om^\star_x$ be the closure of the complement of $I_x(\Om)$: $\Om^\star_x=\overline{(I_x(\Om))^c}$. 
\begin{enumerate}
\item If $x\in\interior\Om$ then $V^{(\a)}_\Om(x)=-V^{(-\a)}_{\Om^\star_x}(x)$ for any $\a$. 
\item If $x\not\in\Om$ then $V^{(\a)}_\Om(x)=V^{(-\a)}_{I_x(\Om)}(x)$ for any $\a$. 
\end{enumerate}
\end{corollary}
When $x=0$, $\Om^\star_0$ is called the {\em star duality} of $\Om$ in \cite{M0}. 
\begin{proof}
Suppose $\ell_{x,v}\ti\partial\Om=\{p_1,\ldots,p_k\}$. Then as 
$$\begin{array}{c}
\displaystyle \ell_{x,v}\ti\partial (I_x(\Om))=\ell_{x,v}\ti\partial\Om^\star_x=\{I_x(p_1),\ldots,I_x(p_k)\}, \\[2mm]
|I_x(p_i)-x|={|p_i-x|}^{-1}, \\[2mm]
-\mbox{\rm sgn}_{I_x(\Om)}(I_x(p_i))=\mbox{\rm sgn}_{\Om^\star_x}(I_x(p_i))=\mbox{\rm sgn}(p_i),
\end{array}$$ 
corollary follows from the theorem above. 
\end{proof}
%
%
\subsection{Log potential}
In the study of $r^\an$-centers in the next section, it is natural to use the log potential 
\[
V_\Om^{\log}(x)=\int_\Om\log\frac1{|x-y|}\,d\mu(y)=-\int_\Om\log{|x-y|}\,d\mu(y)
\]
as $r^0$-potential for the case when $\a=\n$ since 
\[
\frac{\partial V_\Om^{\log}}{\partial x_j}(x)=-\int_\Om{|x-y|}^{-2}(x_j-y_j)\,d\mu(y)
\]
can be considered as 
\[
-\lim_{\a\to\n}\int_\Om{|x-y|}^{\a-\n-2}(x_j-y_j)\,d\mu(y)=-\lim_{\a\to\n}\frac1{\an}\cdot\frac{\partial V_\Om^{(\a)}}{\partial x_j}(x). 
\]

The argument so far in the paper for $V_\Om^{(\a)}$ works roughly as well for $V_\Om^{\log}$. 
It is a continuous function on $\RR^\n$. 
Since 
\[
\textrm{div}_y\left(\frac1\n\left(\log r-\frac1m\right)(y-x)\right)=\log r,
\]
the log potential can also be expressed by the boundary integral as 
\[
V_\Om^{\log}(x)=-\frac1m\int_{\partial\Om}\left(\log r-\frac1m\right)(y-x)\cdot n\,d\sigma(y).
\]
The same argument as in Proposition \ref{prop_derivative_boundary} works to obtain 
\[
\frac{\partial V_\Om^{\log}}{\partial x_j}(x)=-\int_{\partial\Om}\log\frac1{|x-y|}\,\,e_j\cdot n\,d\sigma (y). 
\]
It implies 
\[\begin{array}{rcl}
\displaystyle \frac{\partial^2 V_\Om^{\log}}{\partial {x_j}^2}(x)&=&\displaystyle \int_{\partial\Om}{|x-y|}^{-2}(x_j-y_j)e_j\cdot n\,d\sigma (y),\\[4mm]
\displaystyle \Delta V_\Om^{\log}(x)&=&\displaystyle \int_{\partial\Om}{|x-y|}^{-2}(x-y)\cdot n\,d\sigma (y),
\end{array}
\]
and so, if $\n>2$ 
\[
\Delta V_\Om^{\log}(x)=-(\n-2) V_\Om^{(\n-2)}(x)<0. 
\]
When $\n=2$, as $\Delta_x \log|x-y|=0$ ($|x-y|\ne0$), 
$\Delta V_\Om^{\log}(x)=\Delta V_{D_\e(x)}^{\log}(x)=-2\pi$ if $x\in\interior\Om$ ($\e<\textrm{dist}(x,\partial\Om)$) and $\Delta V_\Om^{\log}(x)=0$ if $x\in\Om^c$.

\section{The $\vect{r^\an}$-centers}
\subsection{Definition and examples}
\begin{definition}\label{def_Mn_center} \rm 
Let $Mm^{(\a)}(\Om)$ denote the minimum value of $V_{\Om}^{(\a)}$ when $\a>\n$, 
the maximum value of $V_\Om^{\log}$ when $\a=\n$, 
the maximum value of $V_{\Om}^{(\a)}$ when $0<\a<\n$, and the maximum value of $V_{\Om}^{(\a)}\big|_{\interior{\Om}}$ when $\a\le0$. 

Define the {\em ${r^{\an}}$-center} of $\Om$, denoted by $C_{\an}(\Om)$, by a point where $Mm^{(\a)}(\Om)$ is attained. 
\end{definition}

Remark that when $\a=\n$ it is meaningless to use $V_{\Om}^{(\n)}$ as it is  constantly equal to $\textrm{Vol}\,(\Om)$. 

\begin{example} \rm 
Let $x=(x_1, \ldots, x_\n)$ be an $r^2$-center of $\Om$. 
Then, as $x$ is a critical point of $V_{\Om}^{(\n+2)}$, it satisfies 
\[0=\frac{\partial }{\partial x_i} \int_{\Om}\, \sum_{k=1}^\n(x_k-y_k)^2 \,dy_1\cdots dy_\n
=2\int_{\Om}\, (x_i-y_i) \,dy_1\cdots dy_\n. \]
Therefore, $x$ is given by 
\[x_i=\frac{\displaystyle \int_{\Om}\, y_i \,dy_1\cdots dy_\n}{\displaystyle \int_{\Om}\, 1 \,dy_1\cdots dy_\n},\] 
which implies that the $r^2$-center coincides with the centroid (center of mass). 
\end{example}

\begin{example} \rm 
In \cite{S} Katsuyuki Shibata introduced the {\em (planar) illuminating center} of a triangle $\Om$ by a point in $\interior\Om$ that gives the maximum value of a map
\[
\Om\setminus(\Om\cap N_{\varepsilon}(\partial\Om))\ni x\mapsto \int_{\Om\setminus B_{\varepsilon}(x)}{{|x-y|}^{-2}}\,d\mu(y)\in\RR
\]
for any sufficiently small $\e>0$, 
where $N_{\varepsilon}(\partial\Om))$ is the $\varepsilon$-tubular neighbourhood of $\partial \Om$. 
In our language, the planar illuminating center is an $r^{-2}$-center. 

Shibata gave the characterization of the planar illuminating center of a triangle as follows (\cite{S}). 
Let $L$ be the illuminating center of a triangle $\triangle ABC$. 
Let $X,Y$ and $Z$ be the intersection points of the edges $BC, CA$ and $AB$ with the line $\overline{AL}, \overline{BL}$, and $\overline{CL}$ respectively (figure \ref{illuminating_center}). 
Then there hold 
\[
\displaystyle \frac{\angle ALB}{\angle ALC} = \frac{|BX|}{|CX|}, \>\, 
\frac{\angle CLA}{\angle CLB} = \frac{|AZ|}{|BZ|}, \,\>\mbox{\rm  and }\>\,
\frac{\angle BLC}{\angle BLA} = \frac{|CY|}{|AY|}.
\]
\begin{figure}[htbp]
\begin{center}
\includegraphics[width=.35\linewidth]{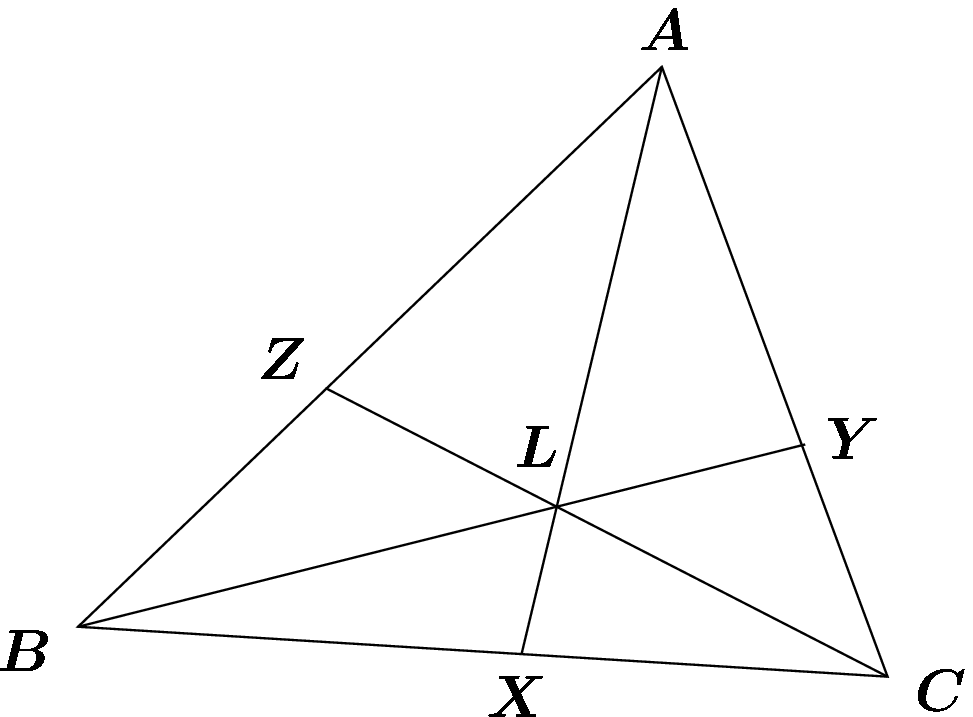}
\caption{Illuminating center of a triangle}
\label{illuminating_center}
\end{center}
\end{figure}
We remark that the last equation of the formulae above can be obtained from the former two by Ceva's theorem. 

We remark that the idea of the illuminating center can be considered as $2$-dimensional analogue of that of the radial center of order $1$ for convex bodies in $\RR^3$ explained in \cite{HMP}. 
\end{example}

\begin{example} \rm 
Let $\Om$ be a convex set with non-empty interior. 
A point $x\in\interior\Om$ is called the {\em radial center of order $\a$} $(\a\ne0)$ if it gives the extreme value of the dual mixed volume of order $\a$ of $\Om_{-x}=\{y-x\,|\,y\in\Om\}$, 
$V_\a(\Om_{-x})=\int_{S^{\n-1}}\big(\rhosub{\Om_{-x}}\big)^\a d\sigma(v)$, 
where $\rho_{\Om_{-x}}$ is the radial function of $\Om_{-x}$. 
It was introduced in \cite{M1} for $0<\a\le1$ and studied in \cite{H} for any $\a$ $(\a\ne0)$. 

The formula \eqref{V^a_int_S^m-1} in Lemma \ref{prop_int_geom} implies that if $\Om$ is convex then an $r^\an$-center $(\a\ne0)$ coincides with the radial center of order $\a$. 
\end{example}

\begin{remark}\rm 
Since $V^{(-\n)}_\Om(x)=-\textrm{\rm Vol}\left(\Om^\star_x\right)$ by Corollary \ref{thm_inversion}, where $\Om^\star_x$ is the closure of the complement of the image of $\Om$ by an inversion in a unit sphere with center $x$, the $r^{-2\n}$-center is a point so that $\Om^\star_x$ has the smallest volume. 
\end{remark}

\subsection{The minimal unfolded regions}\label{subs_min_unfold_region}
\begin{figure}[htbp]
\begin{center}
\includegraphics[width=.32\linewidth]{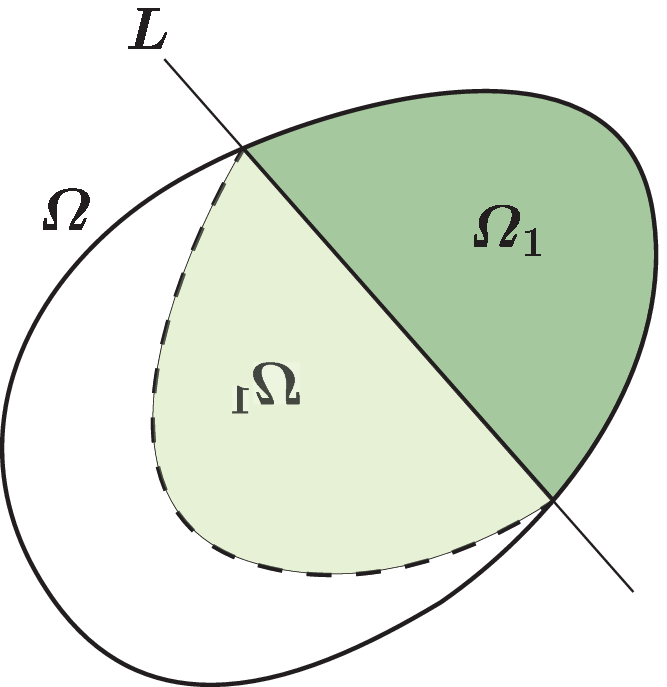}
\caption{Folding a convex set like origami}
\label{origami_convex}
\end{center}
\end{figure}

Let us introduce a region whose complement has no chance to have any ${r^{\an}}$-center. 
Let us first explain in the case when $\Om$ is a convex set in a plan. 
Let $\Om_1$ be a subset of $\Om$ cut out of $\Om$ by a line $L$. 
Fold $\Om$ like origami in $L\cap\Om$. 
Suppose $\Om_1$ can be folded upon $\Om\setminus\Om_1$ (figure \ref{origami_convex}). 
Then our centers cannot be in $\Om_1\setminus(\Om_1\cap L)$ as we will see later. 
%

\begin{definition}\label{def_folded_core} \rm 
Let $\Om$ be a compact set in $\RR^\n$ with piecewise $C^1$ boundary. 
For a unit vector $v$ in $S^{\n-1}$ and a real number $b$, put $M_v=\sup_{x\in\Omega}x\cdot v$, $\Om_{v,b}^+=\Om\cap\{x\in\RR^\n\,|\,x\cdot v> b\}$, and $\textrm{Refl}_{v,b}$ to be a reflection of $\RR^\n$ in a hyperplane $\{x\in\RR^\n\,|\,x\cdot v=b\}$. 
Put 
$$
l_v=\inf\{a\,|\, a\le M_v, \,\textrm{Refl}_{v,b}(\Om_{v,b}^+)\subset\Om\>\>(a\le\forall b\le M_v)\}.
$$
The set $\Om_{v,l_v}^+$ is called the {\em maximal cap} in direction $v$. 
Define the {\em minimal unfolded region} of $\Om$ by 
$$Uf(\Om)
%
%
=\bigcap_{v\in S^{\n-1}}\left\{x\in\RR^\n\,|\,x\cdot v\le l_v\right\}.$$
\end{definition}
%
%
Since the minimal unfolded region contains the centroid, it is a non-empty set. 
As the convex hull of $\Om$ is given by $\bigcap_{\,v\in S^{\n-1}}\{x\,|\,x\cdot v\le M_v\}$, the minimal unfolded region is contained in the convex hull. 
It is compact and convex. 
Remark that the minimal unfolded region of $\Om$ is not necessarily contained in $\Om$ (figure \ref{folded_core_non_convex}). 

\begin{example} \rm 
\begin{enumerate}
\item If $\Om$ is convex and symmetric in a $q$-dimensional hyperplane $H$ ($q<\n$) then the minimal unfolded region is included in $H$. 
Especially, the minimal unfolded region of an $\n$-ball consists of the center. 
\item The minimal unfolded region of a non-obtuse triangle is surrounded by a quadrilateral, two of the edges of which are perpendicular bisectors and the other two are angle bisectors. 
Therefore, two of the vertices of the quadrilateral coincide with the incenter $I$ and the circumcenter $C$ (figure \ref{folded_core_triangle}). 

\begin{figure}[htbp]
\begin{center}
\begin{minipage}{.35\linewidth}
\begin{center}
\includegraphics[width=\linewidth]{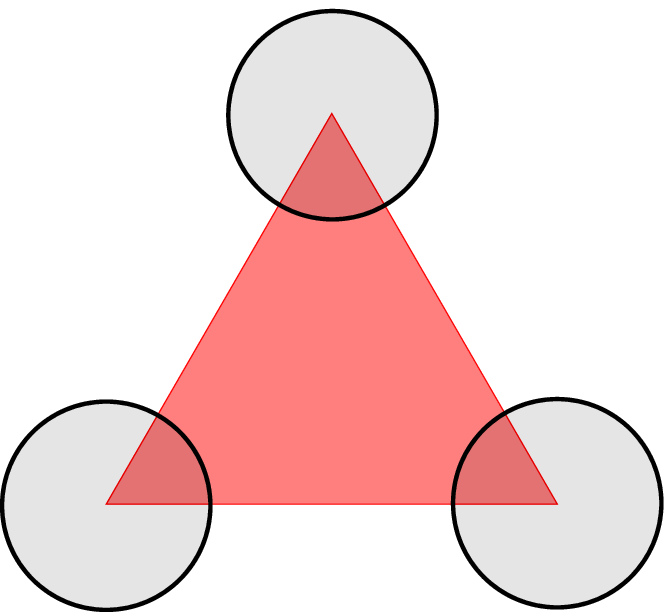}
\caption{A minimal unfolded region (inner triangle) of a non-convex set (union of three discs)}
\label{folded_core_non_convex}
\end{center}
\end{minipage}
\hskip 0.5cm
\begin{minipage}{.5\linewidth}
\begin{center}
\includegraphics[width=\linewidth]{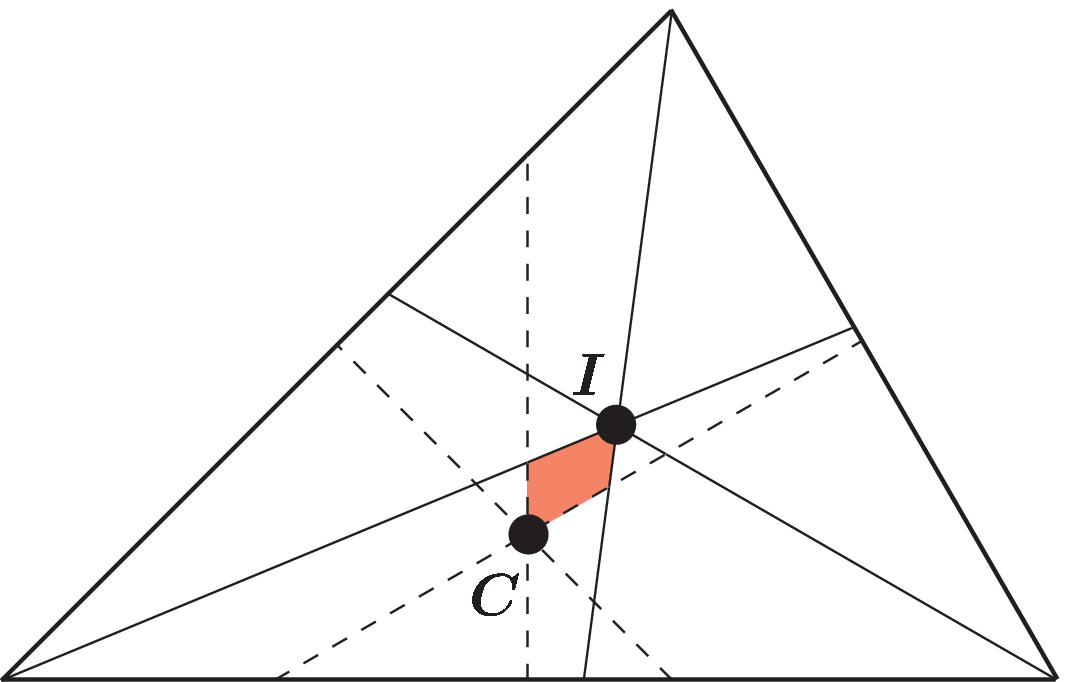}
\caption{A minimal unfolded region of a non-obtuse triangle. Bold lines are angle bisectors, dotted lines are perpendicular bisectors. }
\label{folded_core_triangle}
\end{center}
\end{minipage}
\end{center}
\end{figure}

\end{enumerate}
\end{example}

\begin{problem}\rm 
%
Suppose $\Om$ is convex. 
Is the following true: 
\[
\textrm{diam}(Uf(\Om))
\le\frac12\textrm{diam}(\Om) ?
\]
\end{problem}

\subsection{Existence of $\vect{r^\an}$-centers}
%
\begin{theorem}\label{thm_existence_center_folded_core} 
Any body in $\RR^\n$ with a piecewise $C^1$ boundary that satisfies 
$(\ast\ast)$ in Lemma \ref{lemma_divergence_boundary} has an ${r^{\an}}$-center in the minimal unfolded region for any $\a$. 
\end{theorem}

To be precise, the condition $(\ast\ast)$ is needed for the existence of the center when $\a\le0$. 

\begin{proof} 
(1) First we show that there exists an ${r^{\an}}$-center in the convex hull. 

(i) Suppose $\a>\n$. 
Let $\Om'$ be a convex hull of $\Om$. 
Since $V_{\Om}^{(\a)}$ is continuous on $\RR^\n$ and $\Om'$ is compact, there is a point $x_0$ where the minimum value of $V_{\Om}^{(\a)}$ in $\Om'$ is attained. 
We show that $x_0$ is an $r^\an$-center of $\Om$. 

Let $x$ be a point in $(\Om')^c$. 
Let $x'$ be a point on $\partial\Om'$ so that $\textrm{dist}(x,\partial\Om')=|x-x'|$. 
(The map $x\mapsto x'$ is called the metric projection \cite{M}.) 
Then $\Om'$ is contained in a half-space whose bounary is the hyperplane orthogonal to a line through $x$ and $x'$. 
(This hyperplane is a support hyperplane of $\Om'$ at $x'$.) 
Therefore for any point $y$ in $\Om$ we have $|x-y|>|x'-y|$ and hence $V_{\Om}^{(\a)}(x)>V_{\Om}^{(\a)}(x')\ge V_{\Om}^{(\a)}(x_0)$. 

(ii) The same argument works  for the Riesz potential $(0<\a<\n)$ and the log potential $(\a=\n)$. 

\smallskip
(iii) Suppose $\a\le0$. 
Put $b=\sup_{x\in\interior\Om}V_{\Om}^{(\a)}(x)$. 
As $V_{\Om}^{(\a)}\big|_{\interior\Om}$ is continuous by Proposition \ref{prop_continuity_homothety}, Lemma \ref{lemma_divergence_boundary} (1) implies that 
$\big(V_{\Om}^{(\a)}\big|_{\interior\Om}\big)^{-1}([b-1,b])$ is a compact set in $\RR^\n$. 
Then $b$ is attained there. 

\smallskip
(2) Next we show that an ${r^{\an}}$-center cannot be in the complement of the minimal unfolded region by the moving plane method (\cite{GNN}) using the radial symmetry of partial derivatives of the potential. 
The basic idea is due to Kazuhiro Kurata (personal communication to the author). 
Let us use notation in definition \ref{def_folded_core} in what follows. 

Let $x$ be a critical point of $V_{\Om}^{(\a)}$ when $\a>0$ or $V_{\Om}^{(\a)}\big|_{\interior{\Om}}$ when $\a\le0$. 
Assume $x$ is not in the minimal unfolded region of $\Om$ is 
Then there is a direction $v\in S^{\n-1}$ so that $l_v<x\cdot v\le M_v$. 

Put $c_0=x\cdot v$, $\Om_1=\Om_{v,c_0}^+\cup I_{v,c_0}\left(\Om_{v,c_0}^+\right)$, and $\Om_2=\overline{\Om\setminus\Om_1}$. 
Then the set $\Om$ can be decomposed into two parts, $\Om=\Om_1\cup\Om_2$, 
which is ``disjoint union'' up to the boundaries (figure \ref{folding_unio}). 

\begin{figure}[htbp]
\begin{center}
\includegraphics[width=.4\linewidth]{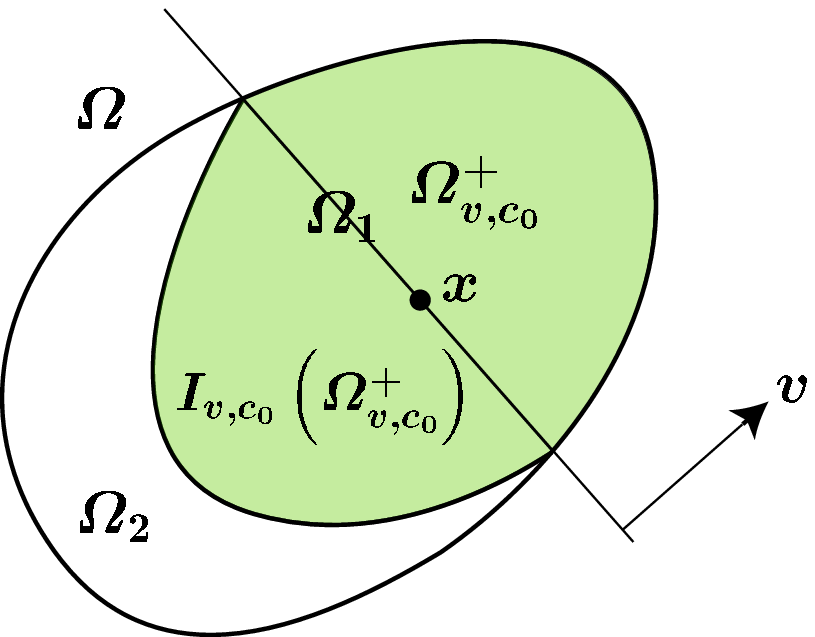}
\caption{}
\label{folding_unio}
\end{center}
\end{figure}

Remark that the slice of $\Om$ with a hyperplane $\{z\in\RR^\n\,|\,z\cdot v=b\}$ grows as $b$ decreases from $M_v$ until $l_v$, namely, $\Om_{v,l_v}^+$ is convex in the direction of $v$. 
Therefore, $\Om_2$ has an non-empty interior as $l_v<c_0\le M_v$. 

Let us first assume that either $x\not\in\partial\Om$ or $\a>1$. 
Then Proposition \ref{prop_derivative_boundary} implies that $V^{(\a)}_\Om$ is partially differentiable at $x$. 
We may assume, after a rotation of $\RR^\n$, that the direction $v$ is equal to $e_1$. 
Then 
\[
\frac{\partial V_{\Om}^{(\a)}}{\partial x_1}(x)
=\frac{\partial V_{\Om_1}^{(\a)}}{\partial x_1}(x)+\frac{\partial V_{\Om_2}^{(\a)}}{\partial x_1}(x). 
\]
The first term of the right hand side vanishes because of the symmetry.  
The second term of the right hand side is given by 
$$
\frac{\partial V_{\Om_2}^{(\a)}}{\partial x_1}(x)
=\int_{\Om_2}(\an)\,{|x-y|}^{\an-2}(x_1-y_1)\,d\mu(y).
$$
The right hand side above is well-defined if $x\not\in\Om_2$ or $\a>1$ since the absolute value of the integrand is bounded above by $|\a-\n|\,{|x-y|}^{\an-1}$. 
Note that if $x\not\in\partial\Om$ then $x\not\in\Om_2$. 

Since $x_1-y_1>0$ on $\interior\Om_2$, which is a non-empty set, $\frac{\partial V_{\Om_2}^{(\a)}}{\partial x_1}(x)$, and hence $\frac{\partial V_{\Om}^{(\a)}}{\partial x_1}(x)$, cannot be equal to $0$, which is a contradiction. 

\smallskip
We are now left with the case when $0<\a\le1$ and $x\in\partial\Om$. 
Take a small neighbourhood $U$ of $x$ which is contained in $\{x'\in\RR^\n\,|\,x'\cdot v>\frac{l_v+c_0}2\}$. 
Then the above argument implies that there is a negative constant $b$ such that if $x'\in U\cap(\partial\Om)^c$ then $\frac{\partial V_{\Om_2}^{(\a)}}{\partial x_1}(x')<b$. 
Therefore, $x$ cannot be a maximum point of $V_{\Om}^{(\a)}$, which is a contradiction. 
\end{proof}
\begin{corollary}\label{cor_center_symmetry} 
\begin{enumerate}
\item The $r^\an$-center of an $\n$-ball is the center for any $\a$. 
\item If $\Om$ is a convex set in a plane and has two lines of symmetry then its $r^\an$-center is the intersection point of the two lines for any $\a$. 
\end{enumerate}
\end{corollary}

\subsection{Extremeness of balls}
We introduce two kinds of extremeness of balls, one applies for any $\a$ whereas the other for $\a>0$, i.e. the case when the renormalization is not needed, but not for any $\a$. 

Recall that $Mm^{(\a)}$ denotes $\min V_{\Om}^{(\a)}(x)$ $(\a>\n)$, $\max V_\Om^{\log}(x)$ $(\a=\n)$, $\max V_{\Om}^{(\a)}(x)$ $0(<\a<\n)$, and $\max V_{\Om}^{(\a)}\big|_{\interior{\Om}}(x)$ $(\a\le0)$. 
%
\begin{proposition}{\rm (Sakata)} 
Suppose $\Om$ has the same volume as a ball $B$. Then we have $Mm^{(\a)}(\Om)\ge Mm^{(\a)}(B)$ if $\a>\n$ and $Mm^{(\a)}(\Om)\le Mm^{(\a)}(B)$ if $\a\le\n$, where the equalities hold if and only if $\Om$ is a ball. 
\end{proposition}
The case when $\n=2$ and $\a>0$ is given in \cite{Sak2}. 
\begin{proof}
The idea is due to \cite{Sak} (see also \cite{GHS}). 

Assume $\a\le0$. Suppose $\Om$ is not a ball. 
By definition, the $r^\an$-center of $\Om$ lies in the interior of $\Om$. 
We may assume that the $r^\an$-centers of $\Om$ and $B$ both coincide with the origin, $\vect 0$. Put $\Om_1=\overline{\Om\setminus(\Om\cap B)}$ and $B_1=\overline{B\setminus (\Om\cap B)}$. %
Lemma \ref{lemma_X-Y} implies that 
$$
Mm^{(\a)}(\Om)-Mm^{(\a)}(B)=
V^{(\a)}_\Om(\vect 0)-V^{(\a)}_B(\vect 0)=V^{(\a)}_{\Om-B}(\vect 0)=
\int_{\Om_1}r^\an d\mu(y)-\int_{B_1}r^\an d\mu(y),
$$
where $r=|y|$. 
As $|y|>|y'|$ for $y\in\interior{\Om_1}$ and $y'\in\interior{B_1}$, $\an<0$, and $\textrm{Vol}(\Om_1)=\textrm{Vol}(B_1)\ne0$, we have 
$$
\int_{\Om_1}r^\an d\mu(y)<\int_{B_1}r^\an d\mu(y),
$$
which implies $Mm^{(\a)}(\Om)< Mm^{(\a)}(B)$. 

The case when $\a>0$ can be proved in the same way. 
\end{proof} 

\smallskip
Define the {\em $r^\an$-energy} of $\Om$ by 
\begin{equation}\label{eq_E}
E^{(\a)}(\Om)=\int_\Om V^{(\a)}_\Om(x)\,d\mu(x)
=\int_{\Om\times\Om} |x-y|^\an\,d\mu(x)d\mu(y)
\end{equation}
for $\a>0$. 

Following Blaschke \cite{B}, T.~Carleman \cite{C} used Steiner's symmetrization to show that if $f(t)$ $(t>0)$ is a positive function with $f'(t)<0$ so that $f(|x-y|)$ is locally integrable in $\RR^\n\times\RR^\n$ then 
$$E_f(\Om)=\int_{\Om\times\Om} f(|x-y|)\,d\mu(x)d\mu(y)$$
attains the maximum only at balls if the volume of $\Om$ is fixed (see also \cite{Mor}). 
%
%
As a corollary it implies 
\begin{proposition}\label{prop_E}
Among bodies with the same volume, balls attain the maximum of $E^{(\a)}(\Om)$ when $0<\a<\n$ and the minimum when $\a>\n$. 
\end{proposition}
The case of gravitational potential ($\n=3, \,\a=2$) was given in \cite{C}. 
\begin{proof}
The first statement can be obtained by putting $f(t)=1/t^{\n-\a}$. 

On the other hand, the second statement can be obtained as follows. 
Suppose $\Om$ has diameter $\diam$. 
Choose a positive function $f(t)$ so that $f(t)=2\diam^\an-t^\an$ for $0<t\le d$ and $f(t)\sim t^{-(\n+1)}$ for $t\ge 2\diam$. Then it follows that $E^{(\a)}(\Om)\ge E^{(\a)}(B)$ 
if $B$ is a ball with the same volume as $\Om$, where the equality holds if only if $\Om$ is a ball. 
\end{proof}
%
%

%
\begin{remark}\rm 
The energy defined by \eqref{eq_E} diverges when $\a\le-1$ since $V^{(\a)}_\Om(x)$ ($x\in\interior\Om$, $\a<0$) can be estimated as $O\left(\left(\textrm{dist}(x,\partial\Om)\right)^\a\right)$ near $\partial\Om$. 
To have a well-defined energy we need one more renormalization process explained in the introduction. Namely, we expand $\displaystyle \int_{\Om\setminus(\Om\cap N_\e(\partial\Om))}V_{\Om}^{(\a)}(x)\,d\mu(x)$ in a series in $\frac1\e$ 
and subtract the divergent terms. 
The detail is studied in \cite{OS} in the case when $\n=2$ and $\a=-2$. 
This renormalized energy can be generalized to Seifert surfaces of knots; thus making a functional on the space of knots, which turns out to be invariant under M\"obius transformations, and can be considered as a renormalization of the average of the square of the linking numbers of the knot and random circles.  

We remark that the proposition above does not hold in general for $\a\le0$. 
For example, when $\n=2$ and $\a=-2$, discs attain the minimum, not the maximum (\cite{OS}). 
\end{remark}

\subsection{Uniqueness of \boldmath $r^\an$-centers}
The $r^\an$-center is not necessarily unique, nor is it always continuous with respect to $\a$. 
Let us see it in a baby case when $\n=1$. 
\begin{example}\label{example_pair_of_intervals} \rm 
Let $\Om$ be a disjoint union of two intervals with the same lengths: $\Om=[-R,-1]\,\cup\,[1,R]$, where $R>1$. 
Theorem \ref{thm_existence_center_folded_core} implies that any $r^{\a-1}$-center belongs to $\big[-\frac{1+R}2,\frac{1+R}2\big]$. 
Observe that if $\a\ge1$ or if $\a<1$ and $x\ne\pm1,\pm R$ then 
$$
\frac{d}{dx}V^{(\a)}_\Om(x)=|x+R|^{\a-1}-|x+1|^{\a-1}-|x-R|^{\a-1}+|x-1|^{\a-1}\,.
$$
\begin{enumerate}
\item Suppose $\a>2$. Then $\frac{d^2}{dx^2}V^{(\a)}_\Om(x)>0$ for any $x$. Therefore the origin is the unique $r^{\a-1}$-center. 
\item Suppose $\a=2$. Then $V^{(2)}_\Om$ is constant on $[-1,1]$ and increases as $|x|$ increases if $|x|>1$. Therefore any point in $[-1,1]$ is an $r^{1}$-center. 
\item Suppose $\a=1$. As $\frac{d}{dx}V^{(1)}_\Om(x)=\log\big|\frac{(R+x)(x-1)}{(R-x)(x+1)}\big|$, which vanishes only when $x=\pm\sqrt{R}$, the $r^{0}$-centers are given by $\pm\sqrt{R}$. 
\item Suppose $\a<2$ $(\a\ne1)$. If $1<\a<2$ then $\frac{d^2}{dx^2}V^{(\a)}_\Om<0$ on $[0,1)$ and $\frac{d^2}{dx^2}V^{(\a)}_\Om>0$ on $(1,R)$. On the other hand, if $\a<1$ then $\frac{d^2}{dx^2}V^{(\a)}_\Om>0$ on $[0,1)$ and $\frac{d^2}{dx^2}V^{(\a)}_\Om<0$ on $(1,R)$. 
In both cases, there is a unique point $x_0=x_0(\a)$ in $(1,R)$ that satisfies $\frac{d}{dx}V^{(\a)}_\Om(x_0(\a))=0$. The $r^{\a-1}$-centers are given by $\pm x_0(\a)$. We remark that Theorem \ref{thm_existence_center_folded_core} implies that $x_0(\a)$ is in $\big(1,\frac{R+1}2\big)$ for any $\a$ ($\a<2, \a\ne1$). 
\end{enumerate}

Let us show that $x_0(\a)$ is a smooth function of $\a$ on $(-\infty,1)\cup(1,2)$. 
Define $f:(1,R)\times\left((-\infty,1)\cup(1,2)\right)\to\RR$ by 
$$
f(x,\a)=\frac{d}{dx}V^{(\a)}_\Om(x)=(x+R)^{\a-1}-(x+1)^{\a-1}-(R-x)^{\a-1}+(x-1)^{\a-1}.
$$
Then 
\[f_x(x,\a)=(\a-1)\left[(x+R)^{\a-2}+(R-x)^{\a-2}-(x+1)^{\a-2}+(x-1)^{\a-2}\right]\ne0\]
since $\a-1\ne0$ and $\displaystyle \frac1{(x-1)^{2-\a}}>\frac1{(x+1)^{2-\a}}$ when $x>1$ and $\a<2$. 
As $f(x,\a)$ is a smooth function of $x$ and $\a$, $x_0(\a)$ is a smooth function of $\a$. 

Note that $\displaystyle \lim_{\a\to2-0}x_0(\a)=1$ as $\displaystyle \lim_{\a\to2-0}f(x,\a)=2x-2$, $\displaystyle \lim_{\a\to1}x_0(\a)=\sqrt{R}$ as 
\[
f(x,\a)=(\a-1)\left(\log\frac{(R+x)(x-1)}{(R-x)(x+1)}+O(\a-1)\right),
\]
and $\displaystyle \lim_{\a\to-\infty}x_0(\a)=\frac{R+1}2$. 
\end{example}

We give a sufficient condition for the uniqueness of the ${r^\an}$-centers. 
\subsubsection{Concavity of \boldmath ${V_\Om^{(\a)}\big|_{\interior\Om}}$ $(\a\le1)$ for convex bodies}
%
\begin{lemma}\label{thm_convexity}
Let $\Om$ be a body in $\RR^\n$ with a piecewise $C^1$ boundary. 
Suppose $\Om$ is convex and $\a\le1$. 
Then $V_{\Om}^{(\a)}$ is a strictly concave function on $\interior\Om$. 
\end{lemma}
\begin{proof}
This is a consequence of Lemma \ref{prop_int_geom} and the following: 
\begin{theorem}{\rm (\cite{HMP}}
If $\phi:\RR_+\to\RR_+$ is concave and strictly increasing, then the map given by 
$$
\RR^\n\ni x\mapsto \int_{S^{\n-1}}\phi(\rhosub{\Om_{-x}}(v))\,d\sigma(v)\in\RR
$$
is strictly concave, where $\rhosub{\Om_{-x}}$ is a radial function of $\Om_{-x}=\{y-x\,|\,y\in\Om\}$. 
\end{theorem}
We can take $\phi(t)=\frac1\a t^\a$ $(\a\le1, \a\ne0)$ or $\phi(t)=\log t$ $(\a=0)$. 
\end{proof}

\begin{theorem}\label{cor_uniqueness_general}
Any convex body in $\RR^\n$ with a piecewise $C^1$ boundary 
has a unique $r^\an$-center in the interior if $\a\le1$. 
\end{theorem}
\begin{proof}
First remark that the conditions $(\ast\ast)$ in Lemma \ref{lemma_divergence_boundary} and $(\ast\ast\ast)$ in Lemma \ref{lemma_partial_differentiability} are satisfied if $\Om$ is convex. 
Therefore, Theorem \ref{thm_existence_center_folded_core} implies the existence of an $r^\an$-center, and Lemma \ref{lemma_partial_differentiability} implies that $r^\an$-center cannot be on $\partial\Om$ when $0<\a\le1$. 
The uniqueness of the maximum point follows from the strong concavity. 
\end{proof}
\begin{remark}\rm 
The statement in \cite{H} that if $\Om$ is convex then $V_{\Om}^{(\a)}\big|_\Om$ (or $V_{\Om}^{(\a)}\big|_{\interior\Om}$ for $\a\le0$) is convex for $\a>\n$ and concave for $\a<\n$ $(\a\ne0)$ does not appear to be correct, as the following counter examples show. 
Nevertheless, we conjecture that Theorem \ref{cor_uniqueness_general} holds for any $\a$. 

Let $\Om$ be an isosceles triangle given by 
$$
\Om=\left\{(x_1,x_2)\,|\,0\le x_1\le 1, 0\le|x_2|\le x_1\tan\left(\frac{\pi}{10}\right)\right\}.
$$

(i) Suppose $\a=1.5$. The graph of the second derivative ${\partial^2 V_{\Om}^{(1.5)}}\big/{\partial x_1{}^2}(t,0)$ produced by Maple using the formula \eqref{f_second_partial_derivative_boundary} is shown in figure \ref{V10_second_deriv}. 
It indicates that $V_{\Om}^{(1.5)}$ is not concave near the vertex $(0,0)$. 
On the other hand, the graph of $V_{\Om}^{(1.5)}(t,0)$ also produced by Maple (figure \ref{V10}) indicates that $V_{\Om}^{(1.5)}$ has a unique maximum point by Corollary \ref{cor_center_symmetry}. 

\begin{figure}[htbp]
\begin{center}
\begin{minipage}{.45\linewidth}
\begin{center}
\includegraphics[width=.8\linewidth]{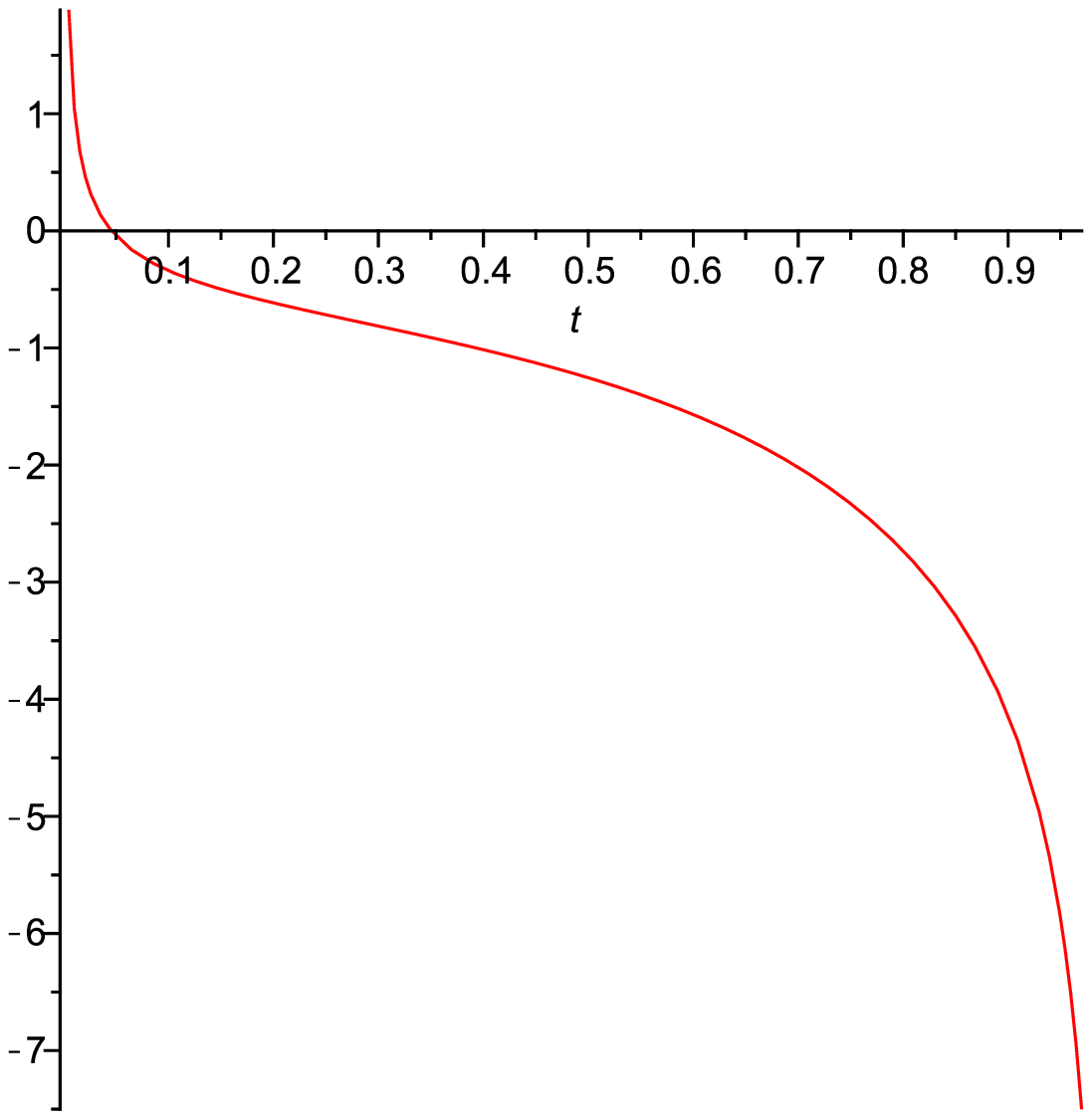}
\caption{$\displaystyle \frac{\partial^2 V_{\Om}^{(1.5)}}{\partial x_1{}^2}\,(t,0)$}
\label{V10_second_deriv}
\end{center}
\end{minipage}
\hskip 0.4cm
\begin{minipage}{.45\linewidth}
\begin{center}
\includegraphics[width=.8\linewidth]{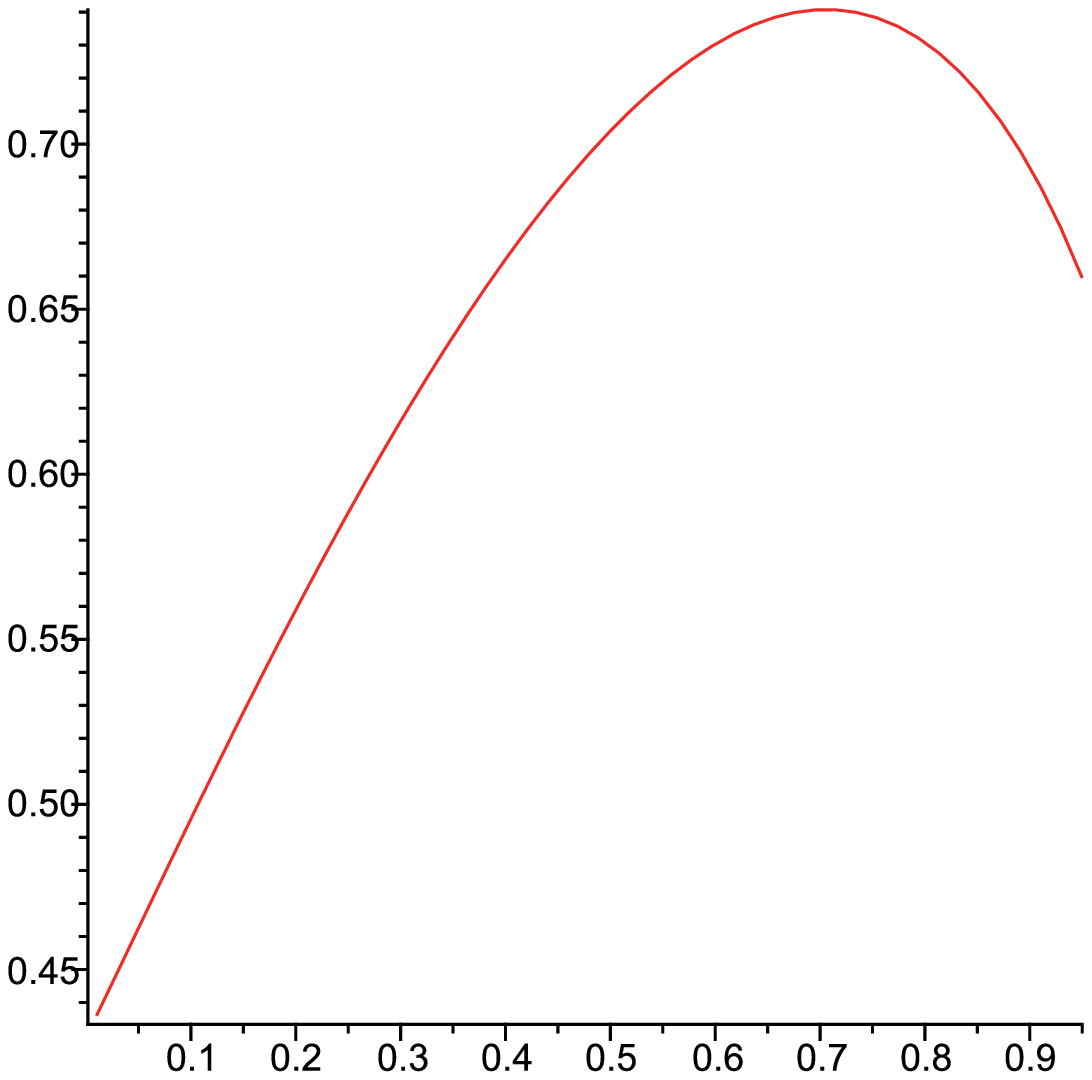}
\caption{$V_{\Om}^{(1.5)}\,(t,0)$}
\label{V10}
\end{center}
\end{minipage}
\end{center}
\end{figure}

(ii) The case when $\a=2.5$ is illustrated in figures \ref{V_2.5_second_deriv} and \ref{V_2.5}. 

\begin{figure}[htbp]
\begin{center}
\begin{minipage}{.45\linewidth}
\begin{center}
\includegraphics[width=.8\linewidth]{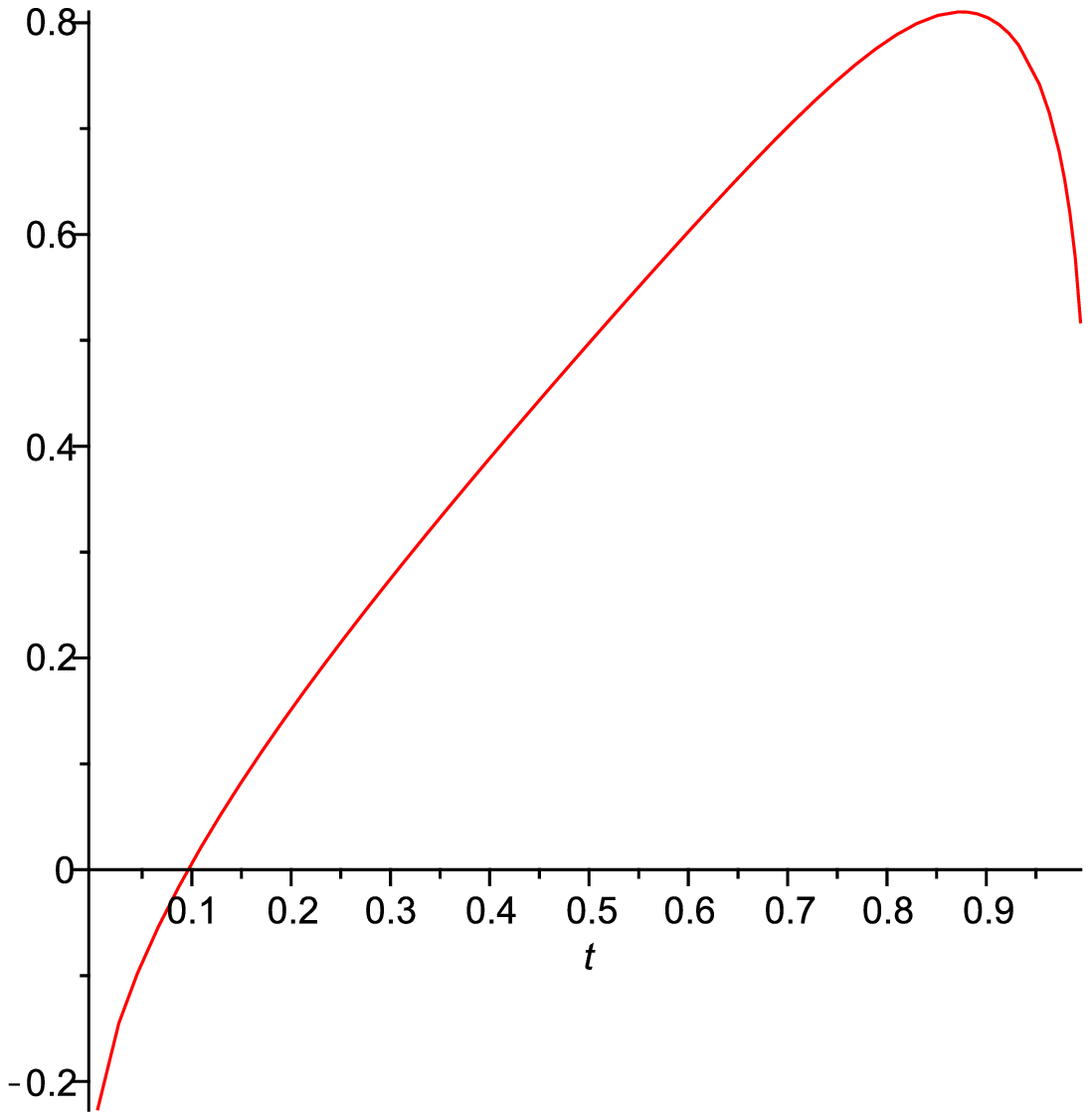}
\caption{$\displaystyle \frac{\partial^2 V_{\Om}^{(2.5)}}{\partial x_1{}^2}\,(t,0)$}
\label{V_2.5_second_deriv}
\end{center}
\end{minipage}
\hskip 0.4cm
\begin{minipage}{.45\linewidth}
\begin{center}
\includegraphics[width=.8\linewidth]{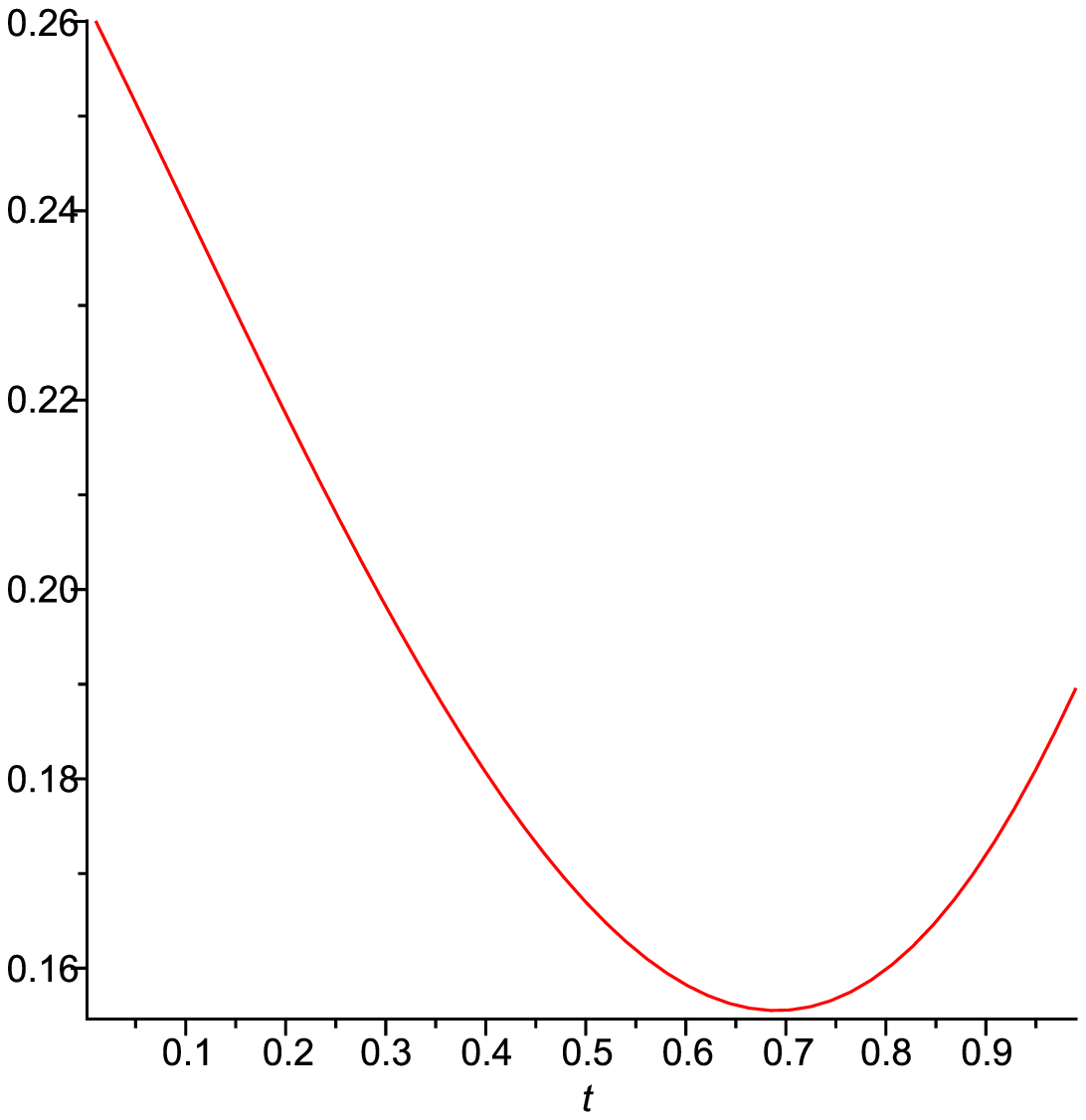}
\caption{$V_{\Om}^{(2.5)}\,(t,0)$}
\label{V_2.5}
\end{center}
\end{minipage}
\end{center}
\end{figure}

\end{remark}

\subsubsection{Convexity of \boldmath $V_\Om^{(\a)}$ $(\a\ge\n+1)$} 
%

\begin{lemma}
Let $\Om$ be a body in $\RR^\n$ with a piecewise $C^1$ boundary. 
If $\a\ge \n+1$ then $\displaystyle \frac{\partial^2 V_{\Om}^{(\a)}}{\partial {x_j}^2}>0$ on $\RR^\n$ for any $j$ $(1\le j\le \n)$. 
\end{lemma}

This follows from the formula \eqref{f_second_partial_derivative_Omega} as $\n\ge2$. 

\begin{theorem}\label{cor_uniquness_alpha_big}
Any body in $\RR^\n$ with a piecewise $C^1$ boundary 
has a unique $r^\an$-center 
if $\a\ge \n+1$. 
\end{theorem}

\begin{proof}
Theorem \ref{thm_existence_center_folded_core} guarantees the existence of a $r^\an$-center. 

Assume there are two $r^\an$-centers $x$ and $x'$. 
We may assume, by a rotation of $\RR^\n$, that the line through $x$ and $x'$ is parallel to the $x_1$-axis. 
As $\frac{\partial V_{\Om}^{(\a)}}{\partial {x_1}}(x)=\frac{\partial V_{\Om}^{(\a)}}{\partial {x_1}}(x')=0$, it contradicts $\frac{\partial^2 V_{\Om}^{(\a)}}{\partial {x_1}^2}>0$. 
\end{proof}

\subsection{Asymptotic behavior of \boldmath $r^\an$-centers as $\a$ goes to $\pm\infty$}
%
In this section we describe the asymptotic behavior of $r^\an$-centers as $\a$ goes to $\pm\infty$. It is generalization of a theorem that an $r^{\a-2}$-center of a non-obtuse triangle approaches the circumcenter as $\a$ goes to $+\infty$ and the incenter as $\a$ goes to $-\infty$, which was announced by Katsuyuki Shibata\footnote{at the 2010 Autumn Meetings of the Mathematical Society of Japan}. 

\begin{lemma}\label{l_asympt_V}
The asymptotic behavior of $V^{(\a)}_\Om(x)$ as $\a$ goes to $\pm\infty$ is given by 
$$
\begin{array}{rcl}
\displaystyle \lim_{\a\to+\infty}\left(V^{(\a)}_\Om(x)\right)^{\frac1\a} &=& \displaystyle \max_{y\in\Om}|y-x|,\\[4mm]
\displaystyle \lim_{\a\to-\infty}\left(-V^{(\a)}_\Om(x)\right)^{-\frac1\a} &=& \displaystyle \frac1{\displaystyle \min_{y\in\overline{{\Om}^c}}|y-x|} \hspace{0.5cm}\big(x\in\interior\Om\big).
\end{array}
$$
\end{lemma}
\begin{proof}
This is because when $\a>0$ we have 
$$
\lim_{\a\to+\infty}\left(V^{(\a)}_\Om(x)\right)^{\frac1\a}=\lim_{\a\to+\infty}\left(\int_\Om|y-x|^\an\,d\mu(y)\right)^{\frac1\a} =  \text{ess}\sup_{\hspace{-0.6cm}y\in\Om}|y-x|=\max_{y\in\Om}|y-x|,
$$
and when $\a<0$ and $x\in\interior\Om$ the formula \eqref{f_Va_complement} implies 
$$
\lim_{\a\to-\infty}\left(-V^{(\a)}_\Om(x)\right)^{-\frac1\a}=\lim_{\a'\to\infty}\left(\int_{\Om^c}\left(\frac1{|y-x|}\right)^{\a'+\n}\,d\mu(y)\right)^{\frac1{\a'}}  =  \text{ess}\sup_{\hspace{-0.6cm}y\in\Om^c}\frac1{|y-x|}=  \frac1{\displaystyle \min_{y\in\overline{{\Om}^c}}|y-x|}.
$$
\end{proof}
\begin{definition} \rm 
(1) We call a point an {\em $r^{\infty}$-center} or a {\em min-max point of $\Om$}, denoted by $C_{\infty}$, if the infimum of the map $\mathbb R^\n\ni x\mapsto \max_{y\in\Om}|y-x|\in\RR$ is attained at the point. In other word, it is the center of a ball with the smallest radius that contains $\Om$. 

\smallskip
(2) We call a point an {\em $r^{-\infty}$-center} or a {\em max-min point of $\Om$}, denoted by $C_{-\infty}$, if the supremum of the map $\mathbb R^\n\ni x\mapsto \min_{y\in\overline{{\Om}^c}}|y-x|\in\RR$ is attained at the point. 
It is a point in $\interior\Om$ where the maximum of the map $\mathbb R^\n\ni x\mapsto \min_{y\in\partial\Om}|y-x|\in\RR$ is attained. In other words, it is a center of a ball with the biggest radius that is contained in $\Om$. 

\smallskip
(3) We call a point an {\em asymptotic $r^{-\infty}$-center of $\Om$} if it is the limit of a convergent sequence of $r^{\a_i-\n}$-centers with $\a_i\to-\infty$ as $i\to+\infty$.  
\end{definition}
It is easy to see that both min-max and max-min points exist. 
Put
$$
R_0
=\min_{x\in\RR^\n}\max_{y\in\Om}|y-x|, 
\>\> r_0 
=\max_{x\in\RR^\n}\min_{y\in\overline{{\Om}^c}}|y-x|. 
$$

\begin{lemma}\label{lemma_unique_min-max}
A min-max point of a compact set $\Om$ is unique, whereas a max-min point is not necessarily unique. 
\end{lemma}
\begin{proof}
The first statement is obvious, and the second can be indicated by a rectangle. 
\end{proof}

\begin{proposition}
A min-max point is contained in the maximal unfolded region, and so is a max-min point if it is unique. 
\end{proposition}
\begin{proof}
Put $R_x=\max_{y\in\Om}|y-x|, \>\>r_x=\min_{y\in\overline{{\Om}^c}}|y-x|$ for $x\in\RR^\n$. 

We use notation in subsection \ref{subs_min_unfold_region}. 
Let $x$ be a point in a maximal cap in direction $v$, $\Om_{v,l_v}^+$. 
Since $\textrm{Refl}_{v,l_v}(\Om_{v,l_v}^+)\subset\Om$, $R_x$ is attained by $y\in\Om\setminus\Om_{v,l_v}^+$, and since $\textrm{Refl}_{v,x\cdot v}(\Om_{v,x\cdot v}^+)\subset\Om$, $r_x$ is attained by $y'\in\partial\Om_{v,x\cdot v}^+\cap\partial\Om$. 
Suppose $x_2=x_1+\e v$ $(x_1,x_2\in\Om_{v,l_v}^+, \e>0)$. 
As an affine hyperplane $\{x\,|\,x\cdot v=l_v\}$ separates $\Om\setminus\Om_{v,l_v}^+$ and $\Om_{v,l_v}^+$, $R_{x_1}<R_{x_2}$, which means that a point in $\Om_{v,l_v}^+$ cannot be a min-max point. 
On the other hand, as $\Om_{v,l_v}^+$ is convex in the direction $v$, $r_{x_1}\ge r_{x_2}$, which means that a point in $\Om_{v,l_v}^+$ cannot be a unique max-min point. 
\end{proof}

\begin{lemma}
Put 
$$
\rhosub{\Om_{-x}}^{\sup}(v)=\sup\{a>0\,|\,x+av\in\Om\}, \>
\rhosub{\Om_{-x}}^{\inf}(v)=\inf\left\{a>0\,|\,x+av\in\overline{\Om^c}\,\right\}
$$
for $x\in\RR^\n,\,v\in S^{\n-1}$ and 
$$
\mathcal{V}_{\rm max}=\{v\in S^{\n-1}\,|\,\rhosub{\Om_{-C_{\infty}}}^{\sup}(v)=R_0\}, \> 
\mathcal{V}_{\rm min}(C_{-\infty})=\{v\in S^{\n-1}\,|\,\rhosub{\Om_{-C_{-\infty}}}^{\inf}(v)=r_0\}, \> 
$$
where $C_{\infty}$ is the unique min-max point and $C_{-\infty}$ a max-min point. 
Then 
$$
\bigcup_{v\in\mathcal{V}_{\rm max}}\{u\in S^{\n-1}\,|\,u\cdot v\ge0\}=
\bigcup_{v\in\mathcal{V}_{\rm min}(C_{-\infty})}\{u\in S^{\n-1}\,|\,u\cdot v\ge0\}=S^{\n-1}. 
$$
\end{lemma}

\begin{proof}
Let us show $\cup_{\stackrel{\phantom{}}{v\in\mathcal{V}_{\rm max}}}\{u\in S^{\n-1}\,|\,u\cdot v\ge0\}=S^{\n-1}$. The proof for $\mathcal{V}_{\min}(C_{-\infty})$ goes parallel. 

Let $B$ be a ball with center $C_{\infty}$ and radius $R_0$, which contains $\Om$. 
If $\partial B$ touches $\partial\Om$ only in points in a compact subset of an open hemisphere of $\partial B$, $\{C_{\infty}+R_0v\,|\,v\in S^{\n-1},\>v\cdot v_0>0\}$ for some $v_0\in S^{\n-1}$, we can move the ball $B$ a little bit in the direction of $v_0$ and make it smaller to have a new ball containing $\Om$ whose radius is smaller than $R_0$, which is a contradiction. 
\end{proof}

It follows that both $\mathcal{V}_{\rm max}$ and $\mathcal{V}_{\rm min}(C_{-\infty})$ contain at least two points each. Therefore, 

\begin{corollary}
Let $\mathcal M$ be a set of max-min points of $\Om$. 
Then it is a subset of the {medial axis}, which is the set of points $x$ in $\interior\Om$ that have at least two points on $\partial\Om$ that give the distance between $x$ and $\partial\Om$. 
\end{corollary}

%
%

\begin{lemma}
For any $\e>0$ there are real numbers $\a_\e^+$ $(\a_\e^+>\n)$ and $\a_\e^-$ $(\a_\e^-<0)$ such that if $\a>\a_\e^+$ then an $r^\an$-center is contained in $B_\e(C_{\infty})$ and that if $\a<\a_\e^-$ then an $r^\an$-center is contained in $N_\e(\mathcal M)$ which is an $\e$-neighbourhood of the set of max-min points $\mathcal M$. 
\end{lemma}

\begin{proof}
We only give the proof of the second statement, as that of the first goes parallel. 

Assume there is no such $\a_\e^-$. 
Then there is a sequence of $r^\an$-centers $\{C_{\a_i-\n}\}_{i\in\mathbb N}$ $(\a_i<0)$ in the complement of $N_\e(\mathcal M)$ with $\a_i\to-\infty$ as $i\to+\infty$. 
As $\{C_{\a_i-\n}\}_{i\in\mathbb N}$ is contained in a compact set $\Om\setminus N_\e(\mathcal M)$, there is a convergent subsequence, which we denote by the same symbol $\{C_{\a_i-\n}\}$, with the limit $C'\in \Om\setminus N_\e(\mathcal M)$. 
We have $r_{C'}<r_0$ as $C'\not\in\mathcal{M}$. 
There is a point $p$ in $\partial\Om$ with $|p-C'|=r_{C'}$. 
Put $\delta=(r_0-r_{C'})/3$ and $b=\textrm{Vol}\left(B_\delta(p)\cap\Om^c\right)$, then $b>0$. 
Suppose $x\in B_\delta(C')\cap\interior\Om$. 
Then if $y\in B_\delta(p)\cap\Om^c$, $|x-y|\le r_{C'}+2\delta=r_0-\delta$, and therefore the formula \eqref{f_Va_complement} implies that $V^{(\a)}_\Om(x)\le-(r_0-\delta)^\an b$ for $\a<0$. 
On the other hand, since $\Om\supset B_{r_0}(C_{-\infty})$, corollary \ref{cor_Omega_1_subset_Omega_2} and the formula \eqref{V^a_int_S^m-1} imply $V^{(\a)}_\Om(C_{-\infty})\ge\frac1\a r_0^\a A(S^{\n-1})$. 
Since 
$$
-\a\left(\frac{r_0}{r_0-\delta}\right)^{-\a}>\frac{(r_0-\delta)^\n A(S^{\n-1})}{b}
$$
for sufficiently small $\a$, $V^{(\a)}_\Om(x)< V^{(\a)}_\Om(C_{-\infty})$ for such $\a$, which implies that a point $x$ in $B_\delta(C')\cap\interior\Om$ cannot be an $r^\an$-center for sufficiently small $\a$. 
This contradicts the assumption that $C'$ is the limit of $r^{\a_i-\n}$-centers with $\a_i\to-\infty$ as $i\to+\infty$. 
\end{proof}

Recall that the {\sl Hausdorff distance} $d_H(X,Y)$ between subsets $X$ and $Y$ in $\RR^\n$ is given by 
$$
d_H(X,Y)=\inf\{\e>0\,|\,X\subset (Y)_\e, \, Y\subset(X)_\e\},\> (X)_\e=\bigcup_{x\in X}B_\e(x).
$$
When $X$ is a singleton, i.e. $X=\{x\}$, $d_H(\{x\},Y)=\sup_{y\in Y}|x-y|$.

\begin{theorem}
\begin{enumerate}
\item The $r^\an$-center converges to the $r^\infty$-center i.e. the min-max point as $\a$ goes to $+\infty$. 

\item An asymptotic $r^{-\infty}$-center is an $r^{-\infty}$-center, i.e. the limit of any convergent sequence of $r^{\a_i-\n}$-centers with $\a_i\to-\infty$ as $i\to+\infty$ is a max-min point. 
Especially, when the max-min point is unique, the sets of $r^\an$-centers $\mathcal{C}_\an$ converges to $\{C_{-\infty}\}$ as $\a$ goes to $-\infty$ with respect to Hausdorff distance. 
\end{enumerate}
\end{theorem}

Department of Mathematics and Informatics, 

Faculty of Science, Chiba University 

1-33 Yayoi-cho, Inage, Chiba, 263-8522, JAPAN. 

E-mail: ohara@math.s.chiba-u.ac.jp


\begin{thebibliography}{E} 
\bibitem[AS]{AS} D.~Auckly, L.~Sadun, { A family of M\"obius invariant 2-knot energies}, Geometric Topology (Athens, GA, 1993), Studies in Advanced Math, AMS, 1997. 

\bibitem[B]{B}W.~Blaschke, { Eine isoperimetrische Eigenschaft des Kreises,} Math. Z. { 1} (1918), 52\,--\,57.


\bibitem[C]{C}T.~Carleman, { \"Uber eine isoperimetrische Aufgabe und ihre physikalischen Anwendungen,} Math. Z. { 3} (1919), 1\,--\,7.

\bibitem[GHS]{GHS}F.~Gao, D.~Hug, R.~Schneider, { Intrinsic volumes and polar sets in spherical space,} Math. Notae { 41} (2003), 159\,--\,176. 

\bibitem[Ga]{Ga}R.J.~Gardner, { Geometric Tomography,} Cambridge Univ. Press, New York, 1995. 

\bibitem[GV]{GV}R.J.~Gardner, A.~Vol\v{c}i\v{c}, { Tomography of convex and star bodies,} Adv. Math. { 108} (1994) 367\,--\,399.

\bibitem[GNN]{GNN}B.~Gidas, W. M. ~Ni, L.~Nirenberg, { Symmetry and related properties via the maximum principle,} Comm. Math. Phys. { 68} (1979), 209\,--\,243.

\bibitem[H]{H}I.~Herburt, { On the Uniqueness of Gravitational Centre,}  Mathematical Physics, Analysis and Geometry { 10} (2007), 251\,--\,259

\bibitem[HMP]{HMP}I.~Herburt, M.~Moszynska, Z.~Peradzynski, { Remarks on radial centres of convex bodies,} Math. Phys. Anal. Geom. { 8} (2005), 157\,--\,172


\bibitem[L1]{Lu75}E.~Lutwak, { Dual mixed volumes,} Pacific J. Math. { 58} (1975), 531\,--\,538. 

\bibitem[L2]{Lu88}E.~Lutwak, { Intersection bodies and dual mixed volumes,} Advances in Mathematics { 71} (1988), 232\,--\,261.

\bibitem[Mor]{Mor}F.~Morgan, { A round ball uniquely minimizes gravitational potential energy,} Proc. Amer. Math. Soc. { 133} (2005), 2733\,--\,2735.

\bibitem[M1]{M0}M.~Moszy\'nska, { . Quotient star bodies, intersection bodies and star duality,} J. Math. Anal. Appl. { 232} (1999), 45\,--\,60.

\bibitem[M2]{M1}M.~Moszy\'nska, { Looking for selectors of star bodies,} Geom. Dedicata { 81} (2000), 131\,--\,147.

\bibitem[M3]{M}M.~Moszy\'nska, { Selected topics in convex geometry,} Birkh\"auser, Boston, 2006.

\bibitem[O1]{O1}J.~O'Hara, { Energy of a knot,} Topology { 30} (1991), 241\,--\,247.

\bibitem[O2]{O2}J.~O'Hara, {  Family of energy functionals of knots,} 
Topology Appl. { 48} (1992), 147\,--\,161.

\bibitem[O3]{O3}J.~O'Hara, { Energy of knots and conformal geometry,} Series on Knots and Everything Vol. 33, World Scientific, Singapore, 2003. 


\bibitem[OS]{OS}J.~O'Hara G.~Solanes, { M\"obius invariant energies and average linking with circles}, arXiv:1010.3764.


\bibitem[Sa]{Sak}S. Sakata, { A volume minimizing problen on the gnomonic projection}, arXiv:1103.4685. 

\bibitem[Sa2]{Sak2}S. Sakata, { On uniqueness of centers of a planar domain with respect to potentials}, arXiv:1103.1727. 


\bibitem[S]{Sch}R.~Schneider, { Convex bodies: the Brunn-Minkowski theory,}  Cambridge University Press, Cambridge, 1993.

\bibitem[Sh]{S}K.~Shibata, { Where should a streetlight be placed in a triangle-shaped park?
Elementary integro-differential geometric optics,} available at 
http://www1.rsp.fukuoka-u.ac.jp/kototoi/shibataaleph-sjs.pdf
\end{thebibliography}
\end{document}